\documentclass[12pt,reqno]{amsart}
\usepackage{cite}
\usepackage{amsfonts}
\usepackage{amssymb, amsmath, latexsym}
\usepackage{longtable}
\usepackage{lipsum}
\usepackage{hyperref}
\usepackage{mleftright}
\usepackage{graphicx,mathtools}
\usepackage{mathtools}

\usepackage{array}

\numberwithin{equation}{section}
\theoremstyle{plain}
\newtheorem{theorem}{Theorem}[section]

\newtheorem{lemma}[theorem]{Lemma}

\newtheorem{corollary}[theorem]{Corollary}

\title[Identities for the Rogers-Ramanujan Continued Fraction] {Identities for the Rogers-Ramanujan Continued Fraction}
\author[N. D. Baruah]{Nayandeep Deka Baruah}
\address{Department of Mathematical Sciences, Tezpur University,  Assam 784028, India}
\email{nayan@tezu.ernet.in}
\author[P. Talukdar]{Pranjal Talukdar}
\address{Department of Mathematical Sciences, Tezpur University,  Assam 784028, India}
\email{pranjaltalukdar113@gmail.com}
\dedicatory{Dedicated to Professor Bruce C. Berndt on the occasion of his 85th birthday}

\allowdisplaybreaks
\begin{document}
\begin{abstract}
We prove some new modular identities for the Rogers\textendash Ramanujan continued fraction. For example, if $R(q)$ denotes the Rogers\textendash Ramanujan continued fraction, then
\begin{align*}&R(q)R(q^4)=\dfrac{R(q^5)+R(q^{20})-R(q^5)R(q^{20})}{1+R(q^{5})+R(q^{20})},\\
&\dfrac{1}{R(q^{2})R(q^{3})}+R(q^{2})R(q^{3})= 1+\dfrac{R(q)}{R(q^{6})}+\dfrac{R(q^{6})}{R(q)},
\end{align*}and\begin{align*}R(q^2)=\dfrac{R(q)R(q^3)}{R(q^6)}\cdot\dfrac{R(q) R^2(q^3) R(q^6)+2  R(q^6) R(q^{12})+ R(q) R(q^3) R^2(q^{12})}{R(q^3) R(q^6)+2  R(q) R^2(q^3)  R(q^{12})+ R^2(q^{12})}.\end{align*}
In the process, we also find some new relations for the Rogers-Ramanujan functions by using dissections of theta functions and the quintuple product identity.
\end{abstract}
\maketitle

\vskip 3mm
\noindent {\footnotesize 2010 Mathematical Reviews Classification Numbers: Primary 11F27, 11P84; Secondary 11A55, 33D90.}

\section{Introduction }\label{Introduction}

For $|q|<1$, the well-known Rogers-Ramanujan continued fraction $R(q)$ is
defined by
$$R(q):=\dfrac{q^{1/5}}{1}_{+}\dfrac{q}{1}_{+}\dfrac{q^2}{1}_{+}\dfrac{q^3}{1}_{+~\cdots,}\quad|q|<1.$$
The Rogers-Ramanujan identities, first proved by Rogers \cite{rogers1894} and then rediscovered by Ramanujan \cite{ramanujan-pcps} are given by  
\begin{align}\label{rrident1}G(q)&:=\sum_{n=0}^\infty\dfrac{q^{n^2}}{(q;q)_n}=\dfrac{1}{(q;q^5)_\infty(q^4;q^5)_\infty}\\\intertext{and} \label{rrident2} H(q)&:=\sum_{n=0}^\infty\dfrac{q^{n(n+1)}}{(q;q)_n}=\dfrac{1}{(q^2;q^5)_\infty(q^3;q^5)_\infty},
\end{align}
where $G(q)$ and $H(q)$ are known as the Rogers-Ramanujan functions. Here and in the sequel, for complex numbers $a$ and $q$ with $|q|<1$, we use the customary notations
\begin{align*}
(a;q)_0&:=1,~(a;q)_n:=\prod_{k=0}^{n-1}(1-aq^k) ~\textup{for}~ n\geq1,~\textup{and}~(a;q)_\infty:=\lim_{n\to \infty}(a;q)_n.
\end{align*}

Rogers \cite{rogers1894} and Ramanujan \cite{nb} (see \cite[Corollary, p. 30]{bcb3}) also proved that  
\begin{align}\label{R-GH}R(q)&=q^{1/5}\dfrac{H(q)}{G(q)}.
\end{align}

For occasional use in the sequel, we set 
\begin{align}
\label{T(q)} T(q):=\dfrac{H(q)}{G(q)}=q^{-1/5}R(q)=\dfrac{(q;q^5)_\infty(q^4;q^5)_\infty}{(q^2;q^5)_\infty(q^3;q^5)_\infty}.
\end{align}

In his first letter to Hardy, written on January 16, 1913, Ramanujan astounded Hardy by an elegant identity writing $R^5(q)$  as a rational expression in $R(q^5)$, namely, 
\begin{align}\label{rrq5}R^5(q)=R(q^5)~\dfrac{1-2R(q^5)+4R^2(q^5)-3R^3(q^5)+R^4(q^5)}{1+R(q^5)+4R^2(q^5)+2R^3(q^5)+R^4(q^5)}.\end{align}
This identity was also recorded by Ramanujan on p. 289 of his second notebook \cite{nb} and p. 365 of his lost notebook \cite{lnb}. To the best of our knowledge, this is the only identity of this sort.
Proofs of \eqref{rrq5} were given by Rogers \cite{rogers1921} in 1921, Watosn \cite{watson-JLMS} in 1929, Ramanathan \cite{ramanathan} in 1984, 
Yi \cite{yi2001} in 2001, and Gugg \cite{gugg2009} in 2009.

In the literature, there are further modular identities relating $R(q)$ with $R(-q)$ and $R(q^n)$ for some positive integers $n$. Ramanujan recorded identities relating $R(q)$ with $R(-q)$, $R(q^2)$, $R(q^3)$ and $R(q^4)$ in the scattered places of his notebooks \cite{nb} and the lost notebook \cite{lnb} and first proofs can be found in \cite{rogers1921}, \cite{andrewsetal}, \cite[Chapter 32]{bcb5}, and \cite[Chapter 1]{andrews-bcb3}. In 1921, Rogers \cite{rogers1921} found a relation connecting 
$R(q)$ with $R(q^{11})$.  Chan and Tan \cite{chan-tan} found a relation connecting $R(q)$ with $R(q^{19})$ in 1999, and Yi \cite{yi2001} found another relation connecting $R(q)$ with $R(q^{7})$ in 2001. We also refer to Trott \cite{trott} for other modular identities discovered ``experimentally" by using Wolfram's \emph{Mathematica}.
 
The main purpose of this paper is to present some new modular identities for $R(q)$. In the process, we also find several new relations for $G(q)$ and $H(q)$. 

In the next theorem and corollary, we present some new identities analogous to \eqref{rrq5}. 
\begin{theorem}\label{theoremnew1}We have
\begin{align}
\label{r(q)r(q^4)} &R(q)R(q^4)=\dfrac{R(q^5)+R(q^{20})-R(q^5)R(q^{20})}{1+R(q^{5})+R(q^{20})},\\
\label{r(q)^2r(q^2) by r(q^4)}&\dfrac{R^2(q)R(q^{2})}{R(q^{4})}=R(q^{5})~\dfrac{R(q^{10})-R(q^{5})R(q^{10})+R(q^{20})}{R(q^{5})R(q^{10})+R(q^{20})+R(q^{5})R(q^{20})},\\
\label{k,r(q)}&\dfrac{1}{R(q)R^2(q^{2})}-{R(q)R^2(q^{2})}=\dfrac{N}{D},
\end{align}
where \begin{align*} D&=R(q^{5})\left(1+R(q^{5})R^2(q^{10})-R^2(q^{5})R^4(q^{10})\right)\\\intertext{and}
N&=1+R(q^{5})+R^2(q^{5})+2R(q^{5})R(q^{10})+2R^2(q^{5})R(q^{10})-2R(q^{5})R^2(q^{10})\\ \notag
&\quad+2R^3(q^{5})R^2(q^{10})-2R^2(q^{5})R^3(q^{10})+2R^3(q^{5})R^3(q^{10})+R^2(q^{5})R^4(q^{10})\\ \notag
&\quad-R^3(q^{5})R^4(q^{10})+R^4(q^{5})R^4(q^{10}).\end{align*}

\end{theorem}

\begin{corollary}\label{cor-new1}
We have
\begin{align}
\label{r(q^4)^2 by r(q)r(q^2)}&\dfrac{R^2(q^4)}{R(q)R(q^{2})}=R(q^{20})~\dfrac{1-R(q^{5})R(q^{10})-R(q^{20})}{R(q^{5})R(q^{10})-R(q^{20})+R(q^{5})R(q^{10})R(q^{20})},\\
\label{r(q)/r(-q)}&\dfrac{R(q)}{R(-q)}=\dfrac{1-R(q^{5})+R(q^{5})R(-q^{5})}{1+R(-q^{5})+R(q^{5})R(-q^{5})},\\
\label{R(q)^2R(q^2)}&\dfrac{R(q^2)}{R(q)^2}-\dfrac{R(q)^2}{R(q^2)} = 4~\dfrac{D}{N},
\end{align}
where $N$ and $D$ are as stated in the previous theorem.
\end{corollary}

%In the following theorem, we recall two of those identities from the lost notebook \cite[p. 205]{lnb} and proved by Andrews and Berndt \cite[pp. 24--25]{andrews-bcb3} and also by Gugg \cite{gugg-ramaja}.
%\begin{theorem}\label{Rqq2q4} We have
%\begin{align}
%\label{RqRq4}&R(q)R(q^4)=\dfrac{R(q^4)-R^2(q)R(q^2)}{R(q^4)+R^2(q^2)},\\
%\label{RqR-qRq4}&R(q)R(-q)R^2(q^2)=\dfrac{R(q)R(-q)-R(q^2)}{R(-q)-R(q)}.
%\end{align}
%\end{theorem}

In the following theorem, we present new addition to the list of modular identities for $R(q)$ available in the literature.

\begin{theorem}\label{theoremnew2} We have
\begin{align}
\label{r(q)1,2,3,6}&\dfrac{1}{R(q^{2})R(q^{3})}+R(q^{2})R(q^{3})= 1+\dfrac{R(q)}{R(q^{6})}+\dfrac{R(q^{6})}{R(q)},\\
\label{r(q)1,2,4,8,16}&R(q^{2})=\dfrac{R^2(q^{4})R(q^{8})-2R(q)R(q^{4})R(q^{16})+R(q^{4})R^2(q^{16})}{R^2(q)R(q^{4})R(q^{8})-2R(q)R(q^{4})R(q^{8})R(q^{16})+R^2(q)R^2(q^{16})},\\
\label{r(q)1,2,3,6a}&R(q^6)= R(q)R^2(q^2)R(q^3)\times\dfrac{R(q)R(q^2)-R(q^3)+R^3(q)}{R^2(q)R(q^3)-R(q^2)R(q^3)+R^3(q)R(q^2)},\\
\label{r(q)1,2,3,6,12}&R(q^2)=\dfrac{R(q)R(q^3)}{R(q^6)}\times\dfrac{R(q) R^2(q^3) R(q^6)+2  R(q^6) R(q^{12})+ R(q) R(q^3) R^2(q^{12})}{R(q^3) R(q^6)+2  R(q) R^2(q^3)  R(q^{12})+ R^2(q^{12})},\\
\label{r(q)1,2,3,4,6,12,24} &R(q^2)=\dfrac{1}{R(q)R(q^3)R(q^6)}\notag\\
&\quad\quad\quad\times\dfrac{2 R(q)R(q^3)R(q^6)R(q^{24})-R(q^4)R(q^6)R^2(q^{12})-R(q^4)R(q^{12})R^2(q^{24})}{R(q)R(q^3)R(q^6)R(q^{12})-2 R(q^4)R^2(q^{12})R(q^{24})+R(q)R(q^3)R^2(q^{24})}.
\end{align}
\end{theorem}

We prove \eqref{r(q)r(q^4)}--\eqref{r(q)1,2,4,8,16} by using dissection formulae for theta functions and certain relations for the Rogers-Ramanujan functions whereas we derive \eqref{r(q)1,2,3,6a}--\eqref{r(q)1,2,3,4,6,12,24}  from some new relations for the Rogers-Ramanujan functions arising from the so-called quintuple product identity (Eq. \eqref{qtpi} in the next section). 

In his brief communication \cite{ramanujan-plms},  Ramanujan gave two algebraic relations between $G(q)$ and $H(q)$, namely, 
\begin{align*}H(q)\{G(q)\}^{11} - q^2G(q)\{H(q)\}^{11} = 1 + 11q\{G(q)H(q)\}^6\\\intertext{and}
H(q)G(q^{11}) - q^2G(q)H(q^{11})=1
\end{align*}
and remarked that each of these formulae is the simplest of a large class. In fact, these two identities are from a set
of forty identities for $G(q)$ and $H(q)$ that Ramanujan never published.  Ramanujan shared some of the identities with Rogers, who proved ten of the identities in \cite{rogers1921}. Watson had Ramanujan's identities and he proved eight identities (two of them already proved by Rogers) in the paper \cite{watson-jims}.  Watson prepared a handwritten list of the forty identities. In 1975, Birch found Watson's handwritten copy of them in the Oxford University Library and published it in \cite{birch}. Watson's handwritten list of Ramanujan's forty identities was finally published in \cite{lnb}. Eventually, the forty identities were proved with a combined effort of Bressoud \cite{bressoud}, Biagioli \cite{biagioli}, Berndt, G. Choi, Y. -S. Choi, Hahn, Yeap, Yee, Yesilyurt and Yi \cite{bcbetal}, and Yesilyurt \cite{yesilyurt1, yesilyurt2}. We refer to \cite[Chapter 8]{andrews-bcb3} for a more comprehensive detail. 

Apart from Ramanujan's forty identities for the Rogers-Ramanujan functions, there are several other identities in the literature. In \cite{bcb-yesilyurt}, Berndt and Yesilyurt established new representations for $G$ and $H$ as linear
combinations of $G$ and $H$ at different arguments and used them in conjunction with some of the previously proved forty identities to prove new identities for the Rogers-Ramanujan functions.  Work on further new identities for $G(q)$ and $H(q)$ can also be found in the papers by Robins \cite{robins}, Chu \cite{chu}, Gugg \cite{gugg2009, gugg-ramaja, gugg-thesis, gugg-jnt, gugg-ijnt}, Koike \cite{koike},  Bringmann and Swisher \cite{bringmann1, bringmann2}, Bulkhali and Dasappa \cite{bulkhali}, and Baruah and Das \cite{ndb-hjd}.

In the following two theorems, we present some new relations for $G(q)$ and $H(q)$. 
\begin{theorem}\label{GH-new}We have
\begin{align}
\label{new 1,2,3,6 first}\dfrac{G\left(q^2\right)G\left(q^3\right)}{H\left(q\right)G\left(q^6\right)}-  \dfrac{H\left(q^2\right)H\left(q^3\right)}{G\left(q\right)H\left(q^6\right)}&= q \dfrac{(q^5;q^5)_\infty^2(q^{30};q^{30})_\infty^2}{(q^{10};q^{10})_\infty^2(q^{15};q^{15})_\infty^2}\\\intertext{and}
\label{new 1,2,3,6 second}\dfrac{H\left(q\right)G\left(q^6\right)}{H\left(q^2\right)H\left(q^3\right)}-  q^2\dfrac{G\left(q\right)H\left(q^6\right)}{G\left(q^2\right)G\left(q^3\right)}&= \dfrac{(q^{10};q^{10})_\infty^2(q^{15};q^{15})_\infty^2}{(q^5;q^5)_\infty^2(q^{30};q^{30})_\infty^2}.
\end{align}
\end{theorem}

\begin{theorem}\label{GH-new-ND-R} Let 
\begin{align*}
S&:=G^2(q)H(q^2)G(q^3)H^2(q^4)H(q^{12}),\\
T&:=H^2(q)G(q^2)H(q^3)G^2(q^4)G(q^{12}),\\
U&:=G^3(q)H^3(q)H(q^6)-H^3(q^2)G(q^3)H(q^3),\\
V&:=G^3(q)H^3(q)G(q^6)-G^3(q^2)G(q^3)H(q^3),\\
 W&:=H^3(q)G(q^6)H(q^6)-G^3(q^2)H^3(q^2)H(q^3),\\
  X&:=G^3(q)G(q^6)H(q^6)-G^3(q^2)H^3(q^2)G(q^3),\\
   Y&:=G^3(q^2)H(q^3)G(q^{12})-H^3(q)G^3(q^4)G(q^6),\\
    Z&:=H^3(q^2)G(q^3)H(q^{12})-G^3(q)H^3(q^4)H(q^6).
\end{align*}
We have
\begin{align}
\label{n1d1}&\dfrac{H(q)H^2(q^3)G(q^{12})+q G(q)G^2(q^3)H(q^{12})}{G(q)G(q^3)H(q^6)G(q^{12})+q^2H(q)H(q^3)G(q^6)H(q^{12})}=\dfrac{\chi(q^5)\chi(-q^{10})}{\chi(-q^{6})\chi(q^{15})},\\
\label{n2d2}&\dfrac{UG(q^2)H(q^3)G(q^{12})-
q VH(q^2)G(q^3)H(q^{12})}{VH(q^2)H(q^6)G(q^{12})
-q^2 UG(q^2)G(q^6)H(q^{12})}=\dfrac{\chi(q^5)\chi(-q^{10})}{\chi(-q^{6})\chi(q^{15})},\\
\label{n3d3}&\dfrac{WG(q)H(q^{12})G(q^{24})-q^2 X H(q)G(q^{12})H(q^{24})
}{q^4WG(q)G(q^6)H(q^{24})
-X H(q)H(q^6)G(q^{24})}=q^3\dfrac{\chi(q^{10})\chi(-q^{12})\chi(-q^{20})\chi(-q^{30})}{\chi^3(-q^{60})},\end{align}
\begin{align}
\label{n4d4}&\dfrac{SYH(q^{12})G(q^{24})-q^2TZG(q^{12})H(q^{24})}{TZH(q^6)G(q^{24})-q^4SYG(q^6)H(q^{24})}=
q^3\dfrac{\chi(q^{10})\chi(-q^{12})\chi(-q^{20})\chi(-q^{30})}{\chi^3(-q^{60})},\\
\label{gh r(q)1,3 by 2,6}&\dfrac{V}{q^{6/5}U}=\dfrac{R(q)R(q^3)}{R^2(q^2)R(q^6)},\\
\label{gh r(q)1,2,3,6}&\dfrac{q^{3/5}W}{X}=R^2(q)R(q^2)R(q^3)R(q^6),\\
\label{gh r(q)1 by 2,4}&\dfrac{Y}{q^3Z}=\dfrac{R(q)}{R(q^2)R(q^4)},
\end{align}
where $\chi(q):=(-q;q^2)_\infty$.
\end{theorem}

We organize the rest of the paper as follows. In Section \ref{sec2}, we present the preliminary results on Ramanujan's theta functions, the Rogers-Ramanujan functions, and some results arising from the quintuple product identity. In Section \ref{sec3}, we prove Theorem \ref{GH-new} by using dissection formulae for theta functions whereas Section \ref{sec4} is devoted to proving Theorem \ref{GH-new-ND-R} by using some theta function identities arising from the quintuple product identity.  In Section \ref{sec5}, we prove Theorem \ref{theoremnew1} and Corollary \ref{cor-new1}. In Section \ref{sec6}, we prove Theorem \ref{theoremnew2}. Partition-theoretic results may be derived from the identities in Theorems \ref{GH-new} and \ref{GH-new-ND-R}. As for examples, in Section \ref{sec7}, we present three such results arising from \eqref{new 1,2,3,6 first}, \eqref{new 1,2,3,6 second}, and \eqref{n1d1}. We offer some concluding remarks in the final section.

\section{Preliminary results and lemmas}\label{sec2}
Ramanujan's general theta function $f(a,b)$ is defined by
\begin{align}
\label{general f(a,b)} f(a,b):=\displaystyle\sum_{n=-\infty}^{\infty}a^{{n(n+1)}/2}b^{{n(n-1)}/2}, \quad \lvert ab\rvert<1.
\end{align}
Jacobi's well-known triple product identity \cite[p. 35, Entry 19]{bcb3} is given by
\begin{align}
\label{JTPI} f(a,b)=(-a;ab)_{\infty}(-b;ab)_{\infty}(ab;ab)_{\infty}.
\end{align}

Three special cases of \eqref{general f(a,b)} are 
\begin{align}\label{phi}\varphi(q)&:=f(q,q)=\sum_{n=-\infty}^\infty q^{n^2}=\dfrac{(q^2;q^2)_\infty^5}{(q;q)_\infty^2(q^4;q^4)_\infty^2},\\
\label{psi}
\psi(q)&:=f(q,q^3)=\sum_{n=0}^\infty q^{n(n+1)/2}=\dfrac{(q^2;q^2)_\infty^2}{(q;q)_\infty},\\
\label{f-q}
f(-q)&:=f(-q,-q^2)=\sum_{n=-\infty}^\infty (-1)^n q^{n(3n-1)/2}=(q;q)_\infty,
\end{align}
where the $q$-product representations arise from \eqref{JTPI} and manipulation of the $q$-products.

In the following two lemmas, we recall some preliminary identities of $f(a,b)$.

\begin{lemma}\textup{(Berndt \cite[pp. 45--46, Entries 29--30]{bcb3})}
\begin{align}\label{f(a,b)+f(-a-b)}f(a,b)&+f(-a,-b)=2f(a^3b,ab^3),\end{align}
\begin{align}\label{f(a,b)-f(-a-b)}f(a,b)&-f(-a,-b)=2af\left(\dfrac{b}{a},a^5b^3\right),\end{align}and if $ab=cd$, then
\begin{align}
\label{f-prod}f(a,b)f(c,d)&=f(ac,bd)f(ad,bc)+af\bigg(\dfrac{b}{c},\dfrac{c}{b}abcd\bigg)f\bigg(\dfrac{b}{d},\dfrac{d}{b}abcd\bigg).
\end{align}
\end{lemma}

\begin{lemma} \textup{(Berndt \cite[p. 51 and p. 350]{bcb3})} We have 
\begin{align}
\label{f(q,q^5)}f(q,q^5)&=\psi(-q^3)\chi(q),\\
\label{f(q,q^2)}f(q,q^2)&=\dfrac{\varphi(-q^3)}{\chi(-q)}.
\end{align}
\end{lemma}

With the aid of the product representations \eqref{rrident1} and \eqref{rrident2} of $G(q)$ and $H(q)$, respectively, and \eqref{JTPI}, the following preliminary identities are apparent (see \cite{bcbetal} and \cite{gugg-jnt} for proofs).
\begin{lemma}\label{GH-prelims}
We have
\begin{align}
&\label{f(-q^2,-q^3)-f(-q,-q^4)}f(-q^2,-q^3)=f(-q)G(q),\quad f(-q,-q^4)=f(-q)H(q),\\
&\label{G(q)H(q)}G(q)H(q)=\frac{f(-q^5)}{f(-q)},\\
\label{f(q^2,q^3)}&f(q^2,q^3)=f(-q^{5}) \frac{H(q)}{H(q^2)}=f(-q) \frac{G(q)H^2(q)}{H(q^2)},\\
\label{f(q,q^4)}&f(q,q^4)=f(-q^{5})\frac{G(q)}{G(q^2)},\\
\label{f(-q,-q^9)}&f(-q,-q^9)=f(-q^{10}) \frac{H(q^2)}{G(q)},\\
\label{f(-q^3,-q^7)}&f(-q^3,-q^7)=f(-q^{10}) \frac{G(q^2)}{H(q)},\\
\label{f(q,q^9)-f(q^3,q^7)}&f(q,q^9)=f(-q^2)G(q)H(q^4),\quad f(q^3,q^7)=f(-q^2)H(q)G(q^4).
\end{align}
\end{lemma}

By manipulating $q$-products, it also follows easily that
\begin{align}
\label{G(-q)-H(-q)}G(-q)&= \dfrac{G(q^2)H^2(q^2)}{G(q)H(q^4)},\quad H(-q)= \dfrac{H(q^2)G^2(q^2)}{H(q)G(q^4)},\\
\label{R(-q)}T(-q)&=\dfrac{T(q^4)}{T(q)T(q^2)},
\end{align}
where $T(q)$ is as defined in \eqref{T(q)}. 

In the following lemma, we recall two useful $5$-dissection formulae from \cite[p. 49 and p. 82]{bcb3}. 
\begin{lemma}We have
\begin{align}
\label{5-dissec phi(q)}\varphi(q)&=\varphi(q^{25})+2q f(q^{15},q^{35})+2q^4f(q^{5},q^{45}),\\
\label{5-dissec psi(q)}\psi(q)&=f(q^{10},q^{15})+qf(q^{5},q^{20})+q^3\psi(q^{25}).
\end{align}
\end{lemma}

The two identities  in the next lemma were recorded by Ramanujan on p. 56 and p. 53, respectively, of his lost notebook \cite{lnb} and proved first by Kang \cite{kang}. 
\begin{lemma}If $u=R(q)$ and $v=R(q^2)$, then
\begin{align}\label{eqn for k,psi(q)}\dfrac{\psi^2(q)}{q\psi^2(q^5)}&=\dfrac{1+uv^2-u^2v^4}{uv^2}\\\intertext{and}
\label{eqn for u,v,psi(q)}\dfrac{\psi(q^{1/5})}{q^{3/5}\psi(q^5)}&=\dfrac{1-uv^2}{uv}+\dfrac{1+uv^2}{v}+1.
\end{align}
\end{lemma}

The identities  in the next lemma were also recorded by Ramanujan on p. 56 of his lost notebook \cite{lnb} and proved first by Kang \cite{kang} (Also see Son \cite{son}). 
\begin{lemma}\label{phiqq5}
We have
\begin{align}
\label{phi(q) plus phi(q^5)}\varphi(q)+\varphi(q^5)&=2q^{4/5}f(q,q^9)R^{-1}(q^4),\\
\vspace{.15cm}
\label{phi(q) minus phi(q^5)}\varphi(q)-\varphi(q^5)&=2q^{1/5}f(q^3,q^7)R(q^4),\\
\vspace{.15cm}
\label{psi(q^2) minus phi(q^{10})} \psi(q^2)-q\psi(q^{10})&=q^{-1/5}f(q^4,q^6)R(q).
\end{align}
\end{lemma}

The identities in the next lemma, which were recorded by Ramanujan in Chapter 19 of his second notebook \cite[Entries 9(iii), 9(vi), 10(iv), 10(v)]{nb}, \cite[p. 258 and p. 262]{bcb3}, readily follow from the previous lemma and \eqref{JTPI}.
\begin{lemma}Let $\chi(q):=(-q;q^2)_{\infty}.$ We have
\begin{align}
\label{phi^2(q) minus phi^2(q^5)} \varphi^2(q)-\varphi^2(q^5)&=4qf(q,q^9)f(q^3,q^7)=4q\chi(q)f(-q^5)f(-q^{20}),\\
\label{psi^2(q) minus q psi^2(q^5)} \psi^2(q)-q\psi^2(q^5)&=f(q,q^4)f(q^2,q^3)=\dfrac{\varphi(-q^5)f(-q^5)}{\chi(-q)}.
\end{align}
\end{lemma}

The following two identities in the next lemma are from the list of forty identities of Ramanujan and they played important roles in proving several other identities in the list; see \cite{watson-jims}, \cite{bcbetal},  \cite{andrews-bcb3} and \cite{yesilyurt2}. The identities can be easily proved by using Lemma \ref{phiqq5} and some preliminary identities of Lemma \ref{GH-prelims}.
\begin{lemma}\label{gh-4} The following identities hold:
 \begin{align}
\label{G(q),G(q^4) plus}G(q)G(q^4)+q H\left(q\right)H\left(q^4\right)&= \dfrac{\varphi(q)}{f(-q^2)},\\
\label{G(q),G(q^4) minus}G\left(q\right)G\left(q^4\right)-q H\left(q\right)H\left(q^4\right)&= \dfrac{\varphi(q^5)}{f(-q^2)}.
\end{align}
\end{lemma}

Next we recall  two more relations on $G(q)$ and $H(q)$.  The first identity is from the list of Ramanujan's forty identities proved first by Rogers \cite{rogers1921} and the second one is from a recent paper by Baruah and Das \cite{ndb-hjd}.
\begin{lemma}We have
\begin{align}
\label{Gq16Hqminus}G(q^{16})H(q)-q^3G(q)H(q^{16})&= \dfrac{\varphi(-q^4)}{f(-q^2)},\\
\label{Gq16Hqplus}G(q^{16})H(q)+q^3G(q)H(q^{16})&= \dfrac{\varphi(-q^{20})}{f(-q^2)}+2q^3\dfrac{f(-q^8)G(q^{16})H(q^{16})}{f(-q^2)}.
\end{align} 
\end{lemma}

In his lost notebook, Ramanujan \cite[p. 207]{lnb}, recorded the quintuple product identity in the form
\begin{align*}
\dfrac{f(-\lambda x,-x^2)}{f(-x,-\lambda x^2)}f(-\lambda x^3)=f(-\lambda^2 x^3,-\lambda x^6)+
x f(-\lambda,-\lambda^2x^9),
\end{align*}
which, by setting $\lambda x^3=q^2$ and $x=-q/B$, may be transformed into (see Berndt \cite[pp. 80--82, Theorem]{bcb3}) 
\begin{align}
\label{qtpi}f(B^3q,q^5/B^3)-B^2f(q/B^3,B^3q^5)=f(-q^2)~\dfrac{f(-B^2,-q^2/B^2)}{f(Bq,q/B)}.
\end{align}

See the excellent paper by Cooper \cite{cooper} for a comprehensive survey of the work on the quintuple product identity and a detailed  analysis of various proofs.

Now we state and prove some theta function identities arising from \eqref{qtpi} which will be used in Section \ref{sec4} to prove Theorem \ref{GH-new-ND-R}.

The following relation on $G(q)$ and $H(q)$ was found by Gugg \cite{gugg-thesis, gugg-jnt} by using some theta function identities arising from \eqref{qtpi}.
\begin{lemma}We have
\begin{align}
\label{G^3(q)H(q^3)}G^3(q)H(q^3)-G(q^3)H^3(q)=\dfrac{3qf^3(-q^{15})}{f(-q)f(-q^3)f(-q^5)}.
\end{align}
\end{lemma}

The  identities in the following lemma were recorded by Ramanujan in his notebooks \cite{nb} and proofs based on \eqref{qtpi} can be found in \cite[p. 379, Entry 10]{bcb3} and \cite[p. 188, Entry 36]{bcb4}.

\begin{lemma}\label{lem-q78}We have
\begin{align}
\label{f(-q^7,-q^8)}f(-q^7,-q^8)+qf(-q^2,-q^{13})&=f(-q^5)~\dfrac{f(-q^2,-q^3)}{f(-q,-q^4)},\\
\label{f(-q^4,-q^{11})}f(-q^4,-q^{11})-qf(-q,-q^{14})&=f(-q^5)~\dfrac{f(-q,-q^4)}{f(-q^2,-q^3)},\\
\label{f(-q^7,-q^8)^3}f^3(-q^7,-q^8)+q^3f^3(-q^2,-q^{13})&=f^3(-q^5)~\dfrac{f(-q^6,-q^9)}{f(-q^3,-q^{12})},\\
\label{f(-q^4,-q^{11})^3} f^3(-q^4,-q^{11})-q^3f^3(-q,-q^{14})&=f^3(-q^5)~\dfrac{f(-q^3,-q^{12})}{f(-q^6,-q^9)}.
\end{align}
\end{lemma}

In the following three lemmas, we state and prove some more useful theta function identities arising from \eqref{qtpi}. 

\begin{lemma}\label{lem-78}We have 
\begin{align}
\label{f(q^7,q^8)}f(q^7,q^8)-qf(q^2,q^{13})&=f(-q^5)~\dfrac{f(-q^2,-q^3)}{f(q,q^4)},\\
\label{f(q^4,q^11)}f(q^4,q^{11})-qf(q,q^{14})&=f(-q^5)~\dfrac{f(-q,-q^4)}{f(q^2,q^3)},\\
\label{f(q^7,q^8)^3}f^3(q^7,q^8)-q^3f^3(q^2,q^{13})&=f^3(-q^5)~\dfrac{f(-q^6,-q^9)}{f(q^3,q^{12})},\\
\label{f(q^4,q^11)^3}f^3(q^4,q^{11})-q^3f^3(q,q^{14})&=f^3(-q^5)~\dfrac{f(-q^3,-q^{12})}{f(q^6,q^9)}.
\end{align}
\end{lemma}

\begin{proof} Setting $q=q^{5/2}$ and $B=q^{3/2}$ in \eqref{qtpi}, we readily arrive at
\begin{align*}
f(q^7,q^8)-qf(q^2,q^{13})=f(-q^5)~\dfrac{f(-q^2,-q^3)}{f(q,q^4)},
\end{align*}
which is \eqref{f(q^7,q^8)}.

Let $\omega=e^{2\pi i/3}$. Putting $q=q^{5/2}$ and $B=\omega q^{3/2}, ~\omega^2 q^{3/2}$, in turn, in \eqref{qtpi}, yields
\begin{align*}
f(q^7,q^8)-q\omega^2f(q^2,q^{13})=f(-q^5)~\dfrac{f(-\omega q^2,-\omega^2q^3)}{f(\omega^2q,\omega q^4)},\\
f(q^7,q^8)-q\omega f(q^2,q^{13})=f(-q^5)~\dfrac{f(-\omega^2 q^2,-\omega q^3)}{f(\omega q,\omega^2 q^4)},
\end{align*}
where the elementary identity $f(a,b)=af\left(a^{-1},a^2b\right)$ has also been used.
Multiplying the previous three identities and then using \eqref{JTPI}, we find that
\begin{align*}
&f^3(q^7,q^8)-q^3f(q^2,q^{13})\\
&=f^3(-q^5)~\dfrac{f(-q^2,-q^3)f(-\omega q^2,-\omega^2q^3)f(-\omega^2 q^2,-\omega q^3)}{f(q,q^4)f(\omega^2q,\omega q^4)f(\omega q,\omega^2 q^4)}\\
&=f^3(-q^5)~\dfrac{(q^2;q^5)_\infty(q^3;q^5)_\infty(\omega^2q^3;q^5)_\infty(\omega q^2;q^5)_\infty(\omega q^3;q^5)_\infty(\omega^2q^2;q^5)_\infty}
{(-q;q^5)_\infty(-q^4;q^5)_\infty(-\omega q^4;q^5)_\infty(-\omega^2 q;q^5)_\infty(-\omega^2 q^4;q^5)_\infty(-\omega q;q^5)_\infty}\\
&=f^3(-q^5)~\dfrac{(q^6;q^{15})_\infty(q^9;q^{15})_\infty}{(-q^3;q^{15})_\infty(-q^{12};q^{15})_\infty}\\
&=f^3(-q^5)~\dfrac{f(-q^6,-q^9)}{f(q^3,q^{12})},
\end{align*}
which is \eqref{f(q^7,q^8)^3}.

Similarly, putting $q=q^{5/2}$ and $B=q^{1/2}$ in \eqref{qtpi}, we arrive at \eqref{f(q^4,q^11)}. Setting $q=q^{5/2}$ and $B=\omega q^{1/2},~ \omega^2 q^{1/2} $ in \eqref{qtpi} and then multiplying the resulting identities and \eqref{f(q^4,q^11)} together, we obtain \eqref{f(q^4,q^11)^3}.
\end{proof}

\begin{lemma}\label{lem-q1317}We have
\begin{align}
\label{f(-q^13,-q^17)}f(-q^{13},-q^{17})+qf(-q^7,-q^{23})&=f(-q^{10})~\dfrac{f(-q^2,-q^8)}{f(-q,-q^9)},\\
\label{f(-q^11,-q^19)}f(-q^{11},-q^{19})+q^3f(-q,-q^{29})&=f(-q^{10})~\dfrac{f(-q^4,-q^6)}{f(-q^3,-q^7)},\\
\label{f(-q^13,-q^17)^3}f^3(-q^{13},-q^{17})+q^3f^3(-q^7,-q^{23})&=f^3(-q^{10})~\dfrac{f(-q^6,-q^{24})}{f(-q^3,-q^{27})},\\
\label{f(-q^11,-q^19)^3}f^3(-q^{11},-q^{19})+q^9f^3(-q,-q^{29})&=f^3(-q^{10})~\dfrac{f(-q^{12},-q^{18})}{f(-q^9,-q^{21})}.
\end{align}
\end{lemma}

\begin{proof}
Setting $q=q^5$ and $B=-q^4$ and $B=-q^2$, in turn, in \eqref{qtpi}, we obtain \eqref{f(-q^13,-q^17)} and \eqref{f(-q^11,-q^19)}.

The identities \eqref{f(-q^13,-q^17)^3} and \eqref{f(-q^11,-q^19)^3} can be proved in a similar way as in the proof of the previous lemma. 
\end{proof}

\begin{lemma}\label{lem-1317}We have 
\begin{align}
\label{f(q^13,q^17)}f(q^{13},q^{17})-qf(q^7,q^{23})&=f(-q^{10})\dfrac{f(-q^2,-q^8)}{f(q,q^9)},\\
\label{f(q^11,q^19)}f(q^{11},q^{19})-q^3f(q,q^{29})&=f(-q^{10})\dfrac{f(-q^4,-q^6)}{f(q^3,q^7)},\\
\label{f(q^13,q^17)^3}f^3(q^{13},q^{17})-q^3f^3(q^7,q^{23})&=f^3(-q^{10})\dfrac{f(-q^6,-q^{24})}{f(q^3,q^{27})},\\
\label{f(q^11,q^19)^3}f^3(q^{11},q^{19})-q^9f^3(q,q^{29})&=f^3(-q^{10})\dfrac{f(-q^{12},-q^{18})}{f(q^9,q^{21})}.
\end{align}
\end{lemma}
\begin{proof}
Putting $q=q^5$ and $B=q^4$ and $B=q^2$, in turn, in \eqref{qtpi} yields \eqref{f(q^13,q^17)} and \eqref{f(q^11,q^19)}. The identities \eqref{f(q^13,q^17)^3} and \eqref{f(q^11,q^19)^3} can be proved by proceeding as in the proof of Lemma \ref{lem-78}.
\end{proof}

\section{Proof of Theorem \ref{GH-new}}\label{sec3} 

\noindent\emph{Proof of \eqref{new 1,2,3,6 first}}.
Setting $a=q$ and $b=q^5$ in \eqref{general f(a,b)}, we have
\begin{align*}
f(q,q^5)=\displaystyle\sum_{n=-\infty}^{\infty}q^{3n^2-2n}.
\end{align*}
It is easily checked that $3n^2-2n \equiv 0, 1$, or $3$ modulo $5$. Therefore, in the series expansion of $f(q,q^5)$, the coefficients of the terms of the forms $q^{5n+2}$ and $q^{5n+4}$ vanish. To exploit this fact further, we try to find a $5$-dissection of $f(q,q^5)$. To that end, first we see from Jacobi triple product identity, \eqref{JTPI}, that
\begin{align}\label{qq5}
f(q,q^5)&=(-q;q^6)_{\infty}(-q^5;q^6)_{\infty}(q^6;q^6)_{\infty}\\
&=\chi(q)\psi(-q^3).\notag
\end{align}

Next, from \eqref{phi^2(q) minus phi^2(q^5)} and \eqref{5-dissec phi(q)}, we have
\begin{align*}
&4q\chi(q)\\
&=\dfrac{1}{f(-q^5)f(-q^{20})}\left(\varphi^2(q)-\varphi^2(q^5)\right)\\
&=\dfrac{1}{f(-q^5)f(-q^{20})}\left(\left(\varphi(q^{25})+2q f(q^{15},q^{35})+2q^4f(q^{5},q^{45})\right)^2-\varphi^2(q^5)\right).
\end{align*}
Employing \eqref{phi^2(q) minus phi^2(q^5)} again in the above, and then dividing both sides by $4q$, we obtain
\begin{align}
\label{proof of new 1,2,3,6 first step1}\chi(q)&=\dfrac{1}{f(-q^5)f(-q^{20})}\bigg(q^4f(q^{5},q^{45})f(q^{15},q^{35})+\varphi(q^{25})f(q^{15},q^{35})\\
&\quad +qf^2(q^{15},q^{35})+q^3\varphi(q^{25})f(q^{5},q^{45})+q^7f^2(q^{5},q^{45})\bigg).\notag
\end{align}

Now, employing the above identity and  \eqref{5-dissec psi(q)} with  $q$ replaced by $-q^3$, in \eqref{qq5}, we obtain the following $5$-dissection of $f(q,q^5)$: 
\begin{align*}
f(q,q^5)&=\dfrac{1}{f(-q^5)f(-q^{20})}\big(q^4f(q^{5},q^{45})f(q^{15},q^{35})+\varphi(q^{25})f(q^{15},q^{35})\\\notag
&\quad +qf^2(q^{15},q^{35})+q^3\varphi(q^{25})f(q^{5},q^{45})+q^7f^2(q^{5},q^{45})\big)\\\notag
&\quad\times \left(f(q^{30},-q^{45})-q^3f(-q^{15},q^{60})-q^9\psi(-q^{75})\right).
\end{align*}
Extracting the terms of the form $q^{5n+2}$ (or $q^{5n+4}$) from both sides of the above after noting from our earlier observation  that no such terms appear on the left-hand side, we find that
\begin{align*}
f(q^{15},q^{35})f(-q^{15},q^{60})
+q^{5}\varphi(q^{25})\psi(-q^{75})-f(q^{5},q^{45})f(q^{30},-q^{45})=0,
\end{align*}
which, by replacing $q^5$ by $-q$, reduces to
\begin{align*}
f(-q^3,-q^7)f(q^3,q^{12})-f(-q,-q^9)f(q^6,q^9)=q\varphi(-q^5)\psi(q^{15}).
\end{align*}
Employing \eqref{phi}, \eqref{psi} and \eqref{f(q^2,q^3)}--\eqref{f(-q^3,-q^7)} in the above, we readily arrive at \eqref{new 1,2,3,6 first}.

\vspace{.3cm}
\noindent\emph{Proof of \eqref{new 1,2,3,6 second}}.
Setting $a=q$ and $b=q^2$ in \eqref{general f(a,b)}, we have
\begin{align*}
f(q,q^2)=\displaystyle\sum_{n=-\infty}^{\infty}q^{(3n^2-n)/2}.
\end{align*}

Observe that $(3n^2-n)/2 \equiv 0, 1$ or $2$ modulo 5.
Thus, in the series expansion of $f(q,q^2)$ the terms of the  forms $q^{5n+3}$ and $q^{5n+4}$ vanish. 

We now find a $5$-dissection of $f(q,q^2)$. By Jacobi triple product identity \eqref{JTPI} and Euler's identity $(-q;q)_\infty=1/(q;q^2)_\infty$, we have
\begin{align}\label{qq2}
f(q,q^2)&=(-q;q^3)_\infty(-q^2;q^3)_\infty(q^3;q^3)_\infty=\dfrac{(-q;q)_\infty(q^3;q^3)_\infty}{(-q^3;q^3)_\infty}\\
&=\dfrac{(q^3;q^3)_\infty^2(q^6;q^6)_\infty}{(q;q^2)_\infty}=\dfrac{\varphi(-q^3)}{\chi(-q)}.\notag
\end{align}

Next, from \eqref{psi^2(q) minus q psi^2(q^5)} and \eqref{5-dissec psi(q)}, we have
\begin{align*}
\dfrac{1}{\chi(-q)}&=\dfrac{1}{f(-q^5)\varphi(-q^5)}\left(\psi^2(q)-q\psi^2(q^5)\right)\\
&=\dfrac{1}{f(-q^5)\varphi(-q^5)}\Big(\left(f(q^{10},q^{15})+qf(q^{5},q^{20})+q^3\psi(q^{25})\right)^2\\
&\quad-q\psi^2(q^5)\Big).
\end{align*}
Once again applying \eqref{psi^2(q) minus q psi^2(q^5)} in the above, we have
\begin{align}
\label{proof of new 1,2,3,6 second step1}\dfrac{1}{\chi(-q)}&=\dfrac{1}{f(-q^5)\varphi(-q^5)}\Big(f^2(q^{10},q^{15})+qf(q^{10},q^{15})f(q^{5},q^{20})+q^2f^2(q^{5},q^{20})\\
&\quad+2q^3\psi(q^{25})f(q^{10},q^{15})+2q^4\psi(q^{25})f(q^{5},q^{20})\Big).\notag
\end{align}

Employing the above identity and \eqref{5-dissec phi(q)}, with  $q$ replaced by $-q^3$, in \eqref{qq2}, we arrive at the following $5$-dissection of $f(q,q^2)$:
\begin{align*}
f(q,q^2)&=\dfrac{1}{f(-q^5)\varphi(-q^5)}\Big(f^2(q^{10},q^{15})+qf(q^{10},q^{15})f(q^{5},q^{20})+q^2f^2(q^{5},q^{20})\\\notag
&\quad+2q^3\psi(q^{25})f(q^{10},q^{15})+2q^4\psi(q^{25})f(q^{5},q^{20})\Big)\\\notag
&\quad\times \Big(\varphi(-q^{75})-2q^3 f(-q^{45},-q^{105})+2q^{12}f(-q^{15},-q^{135})\Big).
\end{align*}

Now, we recall that the series expansion of $f(q,q^2)$ does not contain terms of the forms $q^{5n+3}$ and $q^{5n+4}$. Therefore, extracting the terms of the form $q^{5n+3}$ (or $q^{5n+4}$) from both sides of the above, we find that
\begin{align*}
f(-q^{45},-q^{105})f(q^{10},q^{15})&-q^{10}f(-q^{15},-q^{135})f(q^{5},q^{20})\\
&-\varphi(-q^{75})\psi(q^{25})=0.
\end{align*}
Replacing $q^5$ by $q$ in the above identity, we have
\begin{align*}
f(q^2,q^3)f(-q^9,-q^{21})-q^2f(q,q^4)f(-q^3,-q^{27})=\psi(q^5)\varphi(-q^{15}).
\end{align*}
Employing \eqref{phi}, \eqref{psi}, and \eqref{f(q^2,q^3)}--\eqref{f(-q^3,-q^7)} in the above, we arrive at  \eqref{new 1,2,3,6 second} to finish the proof.

\section{Proof of Theorem \ref{GH-new-ND-R}}\label{sec4} 

\noindent\emph{Proof of \eqref{n1d1}}. 
From \eqref{f(-q^7,-q^8)} and \eqref{f(-q^7,-q^8)^3}, we have
\begin{align}
\label{n1d1 c}&f^3(-q^5)\dfrac{f^3(-q^2,-q^3)}{f^3(-q,-q^4)}-f^3(-q^5)\dfrac{f(-q^6,-q^9)}{f(-q^3,-q^{12})}\\
&=\big(f(-q^7,-q^8)+qf(-q^2,-q^{13})\big)^3-\left(f^3(-q^7,-q^8)+q^3f^3(-q^2,-q^{13})\right)\notag\\
&=3qf(-q^7,-q^8)f(-q^2,-q^{13})\left(f(-q^7,-q^8)+qf(-q^2,-q^{13})\right)\notag\\
&=3qf(-q^7,-q^8)f(-q^2,-q^{13})\bigg(f(-q^5)\dfrac{f(-q^2,-q^3)}{f(-q,-q^4)}\bigg).\notag
\end{align}
 
Applying \eqref{f(-q^2,-q^3)-f(-q,-q^4)} in \eqref{n1d1 c} yields
\begin{align*}
f^2(-q^5)\bigg(\dfrac{G^3(q)}{H^3(q)}-\dfrac{G(q^3)}{H(q^3)}\bigg)=3q\dfrac{G(q)}{H(q)}f(-q^7,-q^8)f(-q^2,-q^{13}).
\end{align*}  
Therefore, 
 \begin{align}
 \label{n1d1 a}G^3(q)H(q^3)-G(q^3)H^3(q)=\dfrac{3q}{f^2(-q^5)}G(q)H^2(q)H(q^3)f(-q^7,-q^8)f(-q^2,-q^{13}).
 \end{align}
 
 In a similar fashion, from \eqref{f(-q^4,-q^{11})}, \eqref{f(-q^4,-q^{11})^3}, and \eqref{f(-q^2,-q^3)-f(-q,-q^4)}, we find that
 \begin{align}
 \label{n1d1 b}G^3(q)H(q^3)-G(q^3)H^3(q)=\dfrac{3q}{f^2(-q^5)}G^2(q)H(q)G(q^3)f(-q^4,-q^{11})f(-q,-q^{14}).
 \end{align}
 
 Comparing \eqref{n1d1 a} and \eqref{n1d1 b}, and then applying \eqref{f-prod}, we arrive at
 \begin{align*}
 &H(q)H(q^3)\left(f(q^9,q^{21})f(q^{10},q^{20})-q^2f(q^6,q^{24})f(q^5,q^{25})\right)\\
 &=G(q)G(q^3)\left(f(q^5,q^{25})f(q^{12},q^{18})-qf(q^3,q^{27})f(q^{10},q^{20})\right).
 \end{align*} 
Employing \eqref{f(q,q^5)}, \eqref{f(q,q^2)}, \eqref{f(q^2,q^3)}, \eqref{f(q,q^4)}, \eqref{f(q,q^9)-f(q^3,q^7)} in the above and then simplifying, we have
\begin{align*}
&\dfrac{f(-q^6)\varphi(-q^{30})}{\chi(-q^{10})}\left(H(q)H^2(q^3)G(q^{12})+qG(q)G^2(q^3)H(q^{12})\right)\\
&=\dfrac{f(-q^{12})f(-q^{30})}{f(-q^{60})}\psi(-q^{15})\chi(q^5)\\
&\quad\times\left(G(q)G(q^3)H(q^6)G(q^{12})+q^2H(q)H(q^3)G(q^6)H(q^{12}) \right),
\end{align*} 
which, with the aid of \eqref{phi}-\eqref{f-q}, can be rearranged to \eqref{n1d1}.

\vspace{.3cm}\noindent\emph{Proof of \eqref{n2d2}}.  
From \eqref{f(q^7,q^8)} and \eqref{f(q^7,q^8)^3}, we have
\begin{align}\label{n2d2 a}
&f^3(-q^5)\left(\dfrac{f(-q^6,-q^9)}{f(q^3,q^{12})}-\dfrac{f^3(-q^2,-q^3)}{f^3(q,q^4)}\right)\\
&=\left(f^3(q^7,q^8)-q^3f^3(q^2,q^{13})\right)-\left(f(q^7,q^8)-qf(q^2,q^{13})\right)^3\notag\\
&=3q f(q^7,q^8)f(q^2,q^{13})\left(f(q^7,q^8)-qf(q^2,q^{13})\right)\notag\\
&=3qf(-q^5) f(q^7,q^8)f(q^2,q^{13})\dfrac{f(-q^2,-q^3)}{f(q,q^4)}.\notag
\end{align}
Applying \eqref{f(-q^2,-q^3)-f(-q,-q^4)} and \eqref{f(q,q^4)} in the above yields
\begin{align}
\label{n2d2 step 1}&G^3(q)H^3(q)G(q^6)-G^3(q^2)G(q^3)H(q^3)\\
&\quad \quad=\dfrac{3}{f^2(-q^5)}G^2(q)H^2(q)G(q^2)G(q^3)H(q^3)f(q^7,q^8)f(q^2,q^{13}).\notag
\end{align}

In a similar way, from \eqref{f(q^4,q^11)} and \eqref{f(q^4,q^11)^3}, we obtain
\begin{align}\label{n2d2 b}
&f^3(-q^5)\left(\dfrac{f(-q^3,-q^{12})}{f(q^6,q^{9})}-\dfrac{f^3(-q,-q^4)}{f^3(q^2,q^3)}\right)=3qf(-q^5) f(q^4,q^{11})f(q,q^{14})\dfrac{f(-q,-q^4)}{f(q^2,q^3)}.
\end{align}

Employing \eqref{f(-q^2,-q^3)-f(-q,-q^4)} and \eqref{f(q^2,q^3)} in \eqref{n2d2 b}, we find that
\begin{align}
\label{n2d2 step 2}&G^3(q)H^3(q)H(q^6)-H^3(q^2)G(q^3)H(q^3)\\
&\quad \quad=\dfrac{3}{f^2(-q^5)}G^2(q)H^2(q)H(q^2)G(q^3)H(q^3)f(q^4,q^{11})f(q,q^{14}).\notag
\end{align}
 
 Dividing \eqref{n2d2 step 1} by \eqref{n2d2 step 2}, we have
 \begin{align*}
 &\dfrac{G^3(q)H^3(q)G(q^6)-G^3(q^2)G(q^3)H(q^3)}{G^3(q)H^3(q)H(q^6)-H^3(q^2)G(q^3)H(q^3)}=\dfrac{G(q^2)f(q^7,q^8)f(q^2,q^{13})}{H(q^2)f(q^4,q^{11})f(q,q^{14})}.
 \end{align*}
With the aid of \eqref{f-prod}, the above identity may be rewritten as
 \begin{align*}
&H(q^2)\bigg(f(q^5,q^{25})f(q^{12},q^{18})+qf(q^{10},q^{20})f(q^{3},q^{27})\bigg)\\
&\quad \quad\times\bigg(G^3(q)H^3(q)G(q^6)-G^3(q^2)G(q^3)H(q^3)\bigg)\\
 &=G(q^2)\bigg(f(q^9,q^{21})f(q^{10},q^{20})+q^2f(q^{6},q^{24})f(q^{5},q^{25})\bigg)\\
 &\quad\quad\times\bigg(G^3(q)H^3(q)H(q^6)-H^3(q^2)G(q^3)H(q^3)\bigg).
 \end{align*}
Employing  \eqref{f(q,q^5)}, \eqref{f(q,q^2)}, \eqref{f(q^2,q^3)}, \eqref{f(q,q^4)}, \eqref{f(q,q^9)-f(q^3,q^7)} in the above and then simplifying further, we arrive at \eqref{n2d2}. 

\vspace{.3cm}
\noindent\emph{Proof of \eqref{n3d3}}.  As the proof is similar to that of  \eqref{n2d2}, we only mention two intermediate identities below.

From \eqref{f(-q^13,-q^17)}--\eqref{f(-q^11,-q^19)^3}, we have
\begin{align*}
&G^3(q)G(q^6)H(q^6)-G^3(q^2)H^3(q^2)G(q^3)\notag\\
& =\dfrac{3q}{f^2(-q^{10})}G(q)G^2(q^2)H^2(q^2)G(q^6)H(q^6)
f(-q^{13},-q^{17})f(-q^7,-q^{23})
\intertext{and}
&H^3(q)G(q^6)H(q^6)-G^3(q^2)H^3(q^2)H(q^3)\notag\\
& =\dfrac{3q^3}{f^2(-q^{10})}H(q)G^2(q^2)H^2(q^2)G(q^6)H(q^6)f(-q^{11},-q^{19})f(-q,-q^{29}).
\end{align*}

\vspace{.3cm}
\noindent\emph{Proof of \eqref{n4d4}}. The proof is also similar to that of \eqref{n2d2}. So we state only two intermediate identities below. 

From \eqref{f(q^13,q^17)}--\eqref{f(q^11,q^19)^3}, it can be shown that
\begin{align*}
&G^3(q)H^3(q^4)H(q^6)-H^3(q^2)G(q^3)H(q^{12})\notag\\
&\quad\quad=\dfrac{3q}{f^2(-q^{10})}G^2(q)H(q^2)G(q^3)H^2(q^4)H(q^{12})
 f(q^{13},q^{17})f(q^7,q^{23})
\intertext{and}
&H^3(q)G^3(q^4)G(q^6)-G^3(q^2)H(q^3)G(q^{12})\notag\\
&\quad\quad=\dfrac{3q^3}{f^2(-q^{10})}H^2(q)G(q^2)H(q^3)G^2(q^4)G(q^{12})
 f(q^{11},q^{19})f(q,q^{29}).
\end{align*}

\vspace{.3cm}\noindent\emph{Proof of \eqref{gh r(q)1,3 by 2,6}}.
From \eqref{n2d2 a}, we have
\begin{align*}&f(-q^6,-q^9)f^3(q,q^4)-f(q^3,q^{12})f^3(-q^2,-q^3)\\
&=\dfrac{3q}{f^2(-q^{5})}f(-q^2,-q^3)f^2(q,q^4)f(q^3,q^{12})f(q^7,q^8)f(q^2,q^{13}).
\end{align*}
Employing the Jacobi triple product identity, \eqref{JTPI}, in the above, we find that
\begin{align}\label{gh r(q)1,3 by 2,6 step 1}
&f(-q^6,-q^9)f^3(q,q^4)-f(q^3,q^{12})f^3(-q^2,-q^3)\\
&=\dfrac{3q}{(q^{5};q^{5})_\infty^2}f(q,q^4)(q^2;q^5)_\infty(q^3;q^5)_\infty(q^5;q^5)_\infty(-q;q^5)_\infty(-q^4;q^5)_\infty(q^5;q^5)_\infty\notag\\
&\quad\times (-q^3;q^{15})_\infty(-q^{12};q^{15})_\infty(q^{15};q^{15})_\infty(-q^7;q^{15})_\infty(-q^8;q^{15})_\infty(q^{15};q^{15})_\infty\notag\\
&\quad\times (-q^2;q^{15})_\infty(-q^{13};q^{15})_\infty(q^{15};q^{15})_\infty\notag\\
&= \dfrac{3qf(-q^2)f^3(-q^{15})}{f(-q)f(-q^{10})}f(q,q^4)f(-q^2,-q^3).\notag
\end{align}
Applying \eqref{f(-q^2,-q^3)-f(-q,-q^4)} and \eqref{f(q,q^4)} in \eqref{gh r(q)1,3 by 2,6 step 1}, we have
\begin{align}
\label{gh r(q)1,3 by 2,6 step 2} \dfrac{G^3(q)H^3(q)G(q^3)}{G^3(q^2)}-\dfrac{G^2(q^3)H(q^3)}{G(q^6)}=\dfrac{3qf(-q^2)f^3(-q^{15})H(q)}{f^2(-q)f(-q^3)f(-q^{10})G(q^2)}.
\end{align}

Next, from \eqref{n2d2 b}, we have
\begin{align*}
&f(-q^3,-q^{12})f^3(q^2,q^3)-f(q^6,q^9)f^3(-q,-q^4)\notag\\
&=\dfrac{3q}{f^2(-q^5)}f(q^4,q^{11})f(q,q^{14})f(-q,-q^4)f^2(q^2,q^3)f(q^6,q^9).
\end{align*}
Applying \eqref{JTPI} in the above, we obtain
\begin{align*}
&f(-q^3,-q^{12})f^3(q^2,q^3)-f(q^6,q^9)f^3(-q,-q^4)\notag\\
&= \dfrac{3qf(-q^2)f^3(-q^{15})}{f(-q)f(-q^{10})}f(-q,-q^4)f(q^2,q^3),
\end{align*}
which by \eqref{f(-q^2,-q^3)-f(-q,-q^4)} and \eqref{f(q^2,q^3)}, reduces to  
\begin{align}
\label{gh r(q)1,3 by 2,6 step 3} \dfrac{G^3(q)H^3(q)H(q^3)}{H^3(q^2)}-\dfrac{G(q^3)H^2(q^3)}{H(q^6)}=\dfrac{3qf(-q^2)f^3(-q^{15})G(q)}{f^2(-q)f(-q^3)f(-q^{10})H(q^2)}.
\end{align}

Dividing \eqref{gh r(q)1,3 by 2,6 step 2} by \eqref{gh r(q)1,3 by 2,6 step 3} and then employing \eqref{R-GH}, we deduce \eqref{gh r(q)1,3 by 2,6}.

\vspace{.3cm}

\noindent\emph{Proof of \eqref{gh r(q)1,2,3,6}}. The proof is similar to that of  \eqref{gh r(q)1,3 by 2,6} given above. Therefore, we mention only two  theta function identities below.

From \eqref{f(-q^13,-q^17)}--\eqref{f(-q^11,-q^19)^3}, we find that
\begin{align*}
&f^3(-q^2,-q^8)f(-q^3,-q^{27})-f^3(-q,-q^9)f(-q^6,-q^{24})\\
&\quad=3q\dfrac{f(-q)f(-q^{10})f^3(-q^{30})}{f(-q^5)}\cdot\dfrac{f(-q,-q^9)}{f(-q^4,-q^6)}\\
\intertext{and}
&f^3(-q^4,-q^6)f(-q^9,-q^{21})-f^3(-q^3,-q^7)f(-q^{12},-q^{18})\\
&\quad=3q^3\dfrac{f(-q)f(-q^{10})f^3(-q^{30})}{f(-q^5)}\cdot\dfrac{f(-q^3,-q^7)}{f(-q^2,-q^8)}.
\end{align*}

\vspace{.3cm}\noindent\emph{Proof of \eqref{gh r(q)1 by 2,4}}. The proof of this identity is also similar to that of  \eqref{gh r(q)1,3 by 2,6} given earlier. So we mention only two important theta function identities.

We obtain the following identities from \eqref{f(q^13,q^17)}--\eqref{f(q^11,q^19)^3} by proceeding as in the proof of \eqref{gh r(q)1,3 by 2,6}.
\begin{align*}
&f^3(q,q^9)f(-q^6,-q^{24})-f^3(-q^2,-q^8)f(q^3,q^{27})\\
&\quad=3q\dfrac{f(-q^2)f(-q^5)f^3(-q^{30})}{f(-q)}\cdot\dfrac{f(-q^2,-q^8)f(q,q^9)}{f(q^2,q^8)f(q^4,q^6)}\\
\intertext{and}
&f^3(q^3,q^7)f(-q^{12},-q^{18})-f^3(-q^4,-q^6)f(q^9,q^{21})\\
&\quad=3q^3\dfrac{f(-q^2)f(-q^5)f^3(-q^{30})}{f(-q)}\cdot\dfrac{f(-q^4,-q^6)f(q^3,q^7)}{f(q^2,q^8)f(q^4,q^6)}.
\end{align*}

\section{Proofs of Theorem \ref{theoremnew1} and Corollary \ref{cor-new1}}\label{sec5} 
\noindent\emph{Proof of \eqref{r(q)r(q^4)}}. 
By \eqref{phi(q) minus phi(q^5)}, we have
\begin{align}
\label{r(q)r(q^4)pf step1}\varphi(q)&=\varphi(q^5)+2q^{1/5}f(q^3,q^7)R(q^4)\\
&=\varphi(q^{25})+2qf(q^{15},q^{35})R(q^{20})+2q^{1/5}f(q^3,q^7)R(q^4).\notag
\end{align}

But, from \eqref{5-dissec phi(q)}, we recall that
\begin{align}
\label{r(q)r(q^4)pf step2}\varphi(q)=\varphi(q^{25})+2q f(q^{15},q^{35})+2q^4f(q^{5},q^{45}).
\end{align}
Comparing \eqref{r(q)r(q^4)pf step1} and \eqref{r(q)r(q^4)pf step2}, we find that
\begin{align*}
q^{1/5}f(q^3,q^7)R(q^4)=  qf(q^{15},q^{35})+q^4f(q^{5},q^{45})-qf(q^{15},q^{35})R(q^{20}).
\end{align*}
Employing \eqref{T(q)} and \eqref{f(q,q^9)-f(q^3,q^7)} in the above, we obtain
\begin{align}
\label{r(q)r(q^4)pf step4} f(-q^2)H(q)H(q^4)=f(-q^{10})\left(H(q^5)G(q^{20})+q^3G(q^5)H(q^{20})-q^4H(q^5)H(q^{20})\right).
\end{align}

Again, by \eqref{phi(q) plus phi(q^5)}, we have
\begin{align}
\label{r(q)r(q^4)pf step5}\varphi(q)&=-\varphi(q^5)+2q^{4/5}f(q,q^9)R^{-1}(q^4)\\&=\varphi(q^{25})-2q^4f(q^5,q^{45})R^{-1}(q^{20})+2q^{4/5}f(q,q^9)R^{-1}(q^4).\notag
\end{align}
Comparing \eqref{r(q)r(q^4)pf step2} and \eqref{r(q)r(q^4)pf step5}, we find that
\begin{align*}
q^{4/5}f(q,q^9)R^{-1}(q^4)=  q^4f(q^5,q^{45})R^{-1}(q^{20})+qf(q^{15},q^{35})+q^4f(q^{5},q^{45}).
\end{align*}
Using \eqref{T(q)} and \eqref{f(q,q^9)-f(q^3,q^7)} in the above, we have
\begin{align*}
&f(-q^2)G(q)G(q^4)\\
&=f(-q^{10})\left(G(q^5)G(q^{20})+qH(q^5)G(q^{20})+q^4G(q^5)H(q^{20})\right).
\end{align*}
It follows from the above identity and \eqref{r(q)r(q^4)pf step4} that
\begin{align*}
\dfrac{H(q)H(q^4)}{G(q)G(q^4)}=\dfrac{H(q^5)G(q^{20})+q^3G(q^5)H(q^{20})-q^4H(q^5)H(q^{20})}{G(q^5)G(q^{20})+qH(q^5)G(q^{20})+q^4G(q^5)H(q^{20})},
\end{align*}
from which, by \eqref{R-GH}, we arrive at
\begin{align*}
R(q)R(q^4)&=\dfrac{R(q^5)+R(q^{20})-R(q^5)R(q^{20})}{1+R(q^{5})+R(q^{20})}.
\end{align*}
Thus, we complete the proof of \eqref{r(q)r(q^4)}.

\vspace{.3cm}
\noindent\emph{Proof of \eqref{r(q)^2r(q^2) by r(q^4)}}.
By \eqref{psi(q^2) minus phi(q^{10})}, we have
\begin{align}
\label{r(q)/r(-q)pf step1}\psi(q^2)&=q\psi(q^{10})+q^{-1/5}f(q^4,q^6)R(q)\\
&=q^6\psi(q^{50})+f(q^{20},q^{30})R(q^5)+q^{-1/5}f(q^4,q^6)R(q).\notag
\end{align}
Comparing \eqref{r(q)/r(-q)pf step1} with \eqref{5-dissec psi(q)}, we find that
\begin{align*}
q^{-1/5}f(q^4,q^6)R(q)=  f(q^{20},q^{30})-f(q^{20},q^{30})R(q^5)+q^2f(q^{10},q^{40}).
\end{align*}
With the aid of \eqref{T(q)}, \eqref{f(q^2,q^3)} and \eqref{f(q,q^4)}, the above identity may be recast as
\begin{align}
\label{r(q)/r(-q)pf step4} &f(-q^2)\dfrac{H(q)G(q^2)H^2(q^2)}{G(q)H(q^4)}\\&=f(-q^{50})\left(\dfrac{H(q^{10})}{H(q^{20})}-q\dfrac{H(q^5)H(q^{10})}{G(q^5)H(q^{20})}+q^2 \dfrac{G(q^{10})}{G(q^{20})}\right).\notag
\end{align}
Replacing $q$ by $-q$ in the above and then employing \eqref{G(-q)-H(-q)}, we find that
\begin{align}
\label{r(q)/r(-q)pf step5}& f(-q^2)\dfrac{G(q)G^2(q^2)H(q^2)}{H(q)G(q^4)}\\&=f(-q^{50})\left(\dfrac{H(q^{10})}{H(q^{20})}+q\dfrac{G(q^5)G(q^{10})}{H(q^5)G(q^{20})}+q^2 \dfrac{G(q^{10})}{G(q^{20})}\right).\notag
\end{align}
Dividing  \eqref{r(q)/r(-q)pf step4} by \eqref{r(q)/r(-q)pf step5}, we find that
\begin{align*}
\dfrac{H^2(q)H(q^2)G(q^4)}{G^2(q)G(q^2)H(q^4)}&=\dfrac{\dfrac{H(q^{10})}{H(q^{20})}-q\dfrac{H(q^5)H(q^{10})}{G(q^5)H(q^{20})}+q^2 \dfrac{G(q^{10})}{G(q^{20})}}{\dfrac{H(q^{10})}{H(q^{20})}+q\dfrac{G(q^5)G(q^{10})}{H(q^5)G(q^{20})}+q^2 \dfrac{G(q^{10})}{G(q^{20})}}.
\end{align*}
By \eqref{R-GH}, the above transforms into 
\begin{align*}
\dfrac{R^2(q)R(q^2)}{R(q^4)}
&=\dfrac{\dfrac{R(q^{10})}{R(q^{20})}-\dfrac{R(q^5)R(q^{10})}{R(q^{20})}+1}{\dfrac{R(q^{10})}{R(q^{20})}+\dfrac{1}{R(q^5)}+1},
\end{align*}
which is clearly equivalent to \eqref{r(q)^2r(q^2) by r(q^4)}.

\vspace{.3cm}
\noindent\emph{Proof of \eqref{k,r(q)}}. Replacing $q$ by $q^5$ in \eqref{eqn for u,v,psi(q)} and then squaring, we have
\begin{align}
\label{k and r(q) pf step 1}\dfrac{\psi^2(q)}{q^6\psi^2(q^{25})}=\left( \dfrac{1-R(q^{5})R^2(q^{10})}{R(q^{5})R(q^{10})}+\dfrac{1+R(q^{5})R^2(q^{10})}{R(q^{10})}+1 \right)^2.
\end{align}

Again, replacing $q$ by $q^5$ in \eqref{eqn for k,psi(q)} and then multiplying the resulting identity with \eqref{eqn for k,psi(q)} again, we find that
\begin{align}
\label{k and r(q) pf step 2}\dfrac{\psi^2(q)}{q^6\psi^2(q^{25})}&=\left(\dfrac{1+R(q^{5})R^2(q^{10})-R^2(q^{5})R^4(q^{10})}{R(q^{5})R^2(q^{10})}\right)\left(1+\dfrac{1}{uv^2}-{uv^2}\right),
\end{align}
where $u=R(q)$ and $v=R(q^2)$. From \eqref{k and r(q) pf step 1} and \eqref{k and r(q) pf step 2}, it follows that
\begin{align*}
\dfrac{1}{uv^2}-{uv^2}&=\left(\dfrac{R(q^{5})R^2(q^{10})}{1+R(q^{5})R^2(q^{10})-R^2(q^{5})R^4(q^{10})}\right)\\\notag&\quad\times\left( \dfrac{1-R(q^{5})R^2(q^{10})}{R(q^{5})R(q^{10})}+\dfrac{1+R(q^{5})R^2(q^{10})}{R(q^{10})}+1 \right)^2-1,
\end{align*}
which upon simplification gives \eqref{k,r(q)}. This completes the proof of Theorem \ref{theoremnew1}.

\vspace{.3cm}
Now we prove the identities in Corollary \ref{cor-new1}.

\vspace{.3cm}
\noindent\emph{Proofs of \eqref{r(q^4)^2 by r(q)r(q^2)}} and \eqref{r(q)/r(-q)}.
Replacing $q$ by $-q$ in \eqref{r(q)r(q^4)} and then applying \eqref{R(-q)} and \eqref{T(q)}, we easily deduce \eqref{r(q^4)^2 by r(q)r(q^2)}. Similarly, one can deduce \eqref{r(q)/r(-q)} by using \eqref{R(-q)} and \eqref{T(q)} in \eqref{r(q)^2r(q^2) by r(q^4)}.

\vspace{.3cm}
\noindent\emph{Proof of \eqref{R(q)^2R(q^2)}}. On p. 326 in his second notebook \cite{nb}, \cite[p. 12]{bcb5}, Ramanujan recorded the following modular equation relating $R(q)$ and $R(q^2)$:
\begin{align*}
R(q)R(q^2)^2= \dfrac{R(q^2)-R^2(q)}{R(q^2)+R^2(q)}.
\end{align*}
This identity was proved first by Rogers \cite{rogers1921}. 
%As shown by Gugg \cite{gugg-ramaja}, the identity also follows easily from \eqref{robins1} and \eqref{robins2}. 
 
The above identity may be recast in the form
\begin{align*}
\dfrac{1}{R(q)R(q^2)^2}- R(q)R(q^2)^2=\dfrac{4}{\dfrac{R(q^2)}{R(q)^2}-\dfrac{R(q)^2}{R(q^2)}}.
\end{align*}
From the above identity and \eqref{k,r(q)}, we readily deduce  \eqref{R(q)^2R(q^2)} to finish the proof of Corollary \ref{cor-new1}.

\section{Proof of Theorem \ref{theoremnew2}}\label{sec6}

\noindent\emph{Proof of \eqref{r(q)1,2,3,6}}.
Multiplying \eqref{new 1,2,3,6 first}
 and \eqref{new 1,2,3,6 second}, we have
 \begin{align*}
 \dfrac{G(q^2)G(q^3)}{H(q^2)H(q^3)}+q^2\dfrac{H(q^2)H(q^3)}{G(q^2)G(q^3)}-\dfrac{G(q^6)H(q)}{G(q)H(q^6)}-q^2\dfrac{G(q)H(q^6)}{G(q^6)H(q)}=q,
 \end{align*}
which, by \eqref{R-GH}, readily transforms into \eqref{r(q)1,2,3,6}.

\vspace{.3cm}
\noindent\emph{Proof of \eqref{r(q)1,2,4,8,16}}. Dividing \eqref{Gq16Hqplus}  by \eqref{Gq16Hqminus}, we have
\begin{align*}
\dfrac{G(q^{16})H(q)+q^3G(q)H(q^{16})}{G(q^{16})H(q)-q^3G(q)H(q^{16})}=\dfrac{\varphi(-q^{20})+2q^3f(-q^8)G(q^{16})H(q^{16})}{\varphi(-q^{4})}.
\end{align*}
Replacing $q$ by $-q^{4}$ in \eqref{G(q),G(q^4) plus} and \eqref{G(q),G(q^4) minus}, and then using the resulting identities in the above identity, we find that
\begin{align*}
&\dfrac{G(q^{16})H(q)+q^3G(q)H(q^{16})}{G(q^{16})H(q)-q^3G(q)H(q^{16})}\\
&=\dfrac{G(-q^4)G(q^{16})+q^4H(-q^4)H(q^{16})+2q^3G(q^{16})H(q^{16})}{G(-q^4)G(q^{16})-q^4H(-q^4)H(q^{16})},
\end{align*}
which may be rewritten with the aid of  \eqref{T(q)} as
\begin{align}\label{r(q)1,2,4,8,16 step1}
&\dfrac{T(q)+q^3T(q^{16})}{T(q)-q^3T(q^{16})}=\dfrac{G(-q^4)G(q^{16})+q^4H(-q^4)H(q^{16})+2q^3G(q^{16})H(q^{16})}{G(-q^4)G(q^{16})-q^4H(-q^4)H(q^{16})}.
\end{align}
Replacing $q$ by $-q$ in \eqref{r(q)1,2,4,8,16 step1}, we have
\begin{align}
\label{r(q)1,2,4,8,16 step2}
\dfrac{T(-q)-q^3T(q^{16})}{T(-q)+q^3T(q^{16})}=\dfrac{G(-q^4)G(q^{16})+q^4H(-q^4)H(q^{16})-2q^3G(q^{16})H(q^{16})}{G(-q^4)G(q^{16})-q^4H(-q^4)H(q^{16})}.
\end{align}
Adding \eqref{r(q)1,2,4,8,16 step1} and \eqref{r(q)1,2,4,8,16 step2}, and then using \eqref{T(q)}, we find that
\begin{align*}
\dfrac{T(q)+q^3T(q^{16})}{T(q)-q^3T(q^{16})}+\dfrac{T(-q)-q^3T(q^{16})}{T(-q)+q^3T(q^{16})}&=2\dfrac{G(-q^4)G(q^{16})+q^4H(-q^4)H(q^{16})}{G(-q^4)G(q^{16})-q^4H(-q^4)H(q^{16})}\\\notag
&=2\dfrac{1+q^4T(-q^4)T(q^{16})}{1-q^4T(-q^4)T(q^{16})}.
\end{align*}
Employing \eqref{R(-q)} in the above and then simplifying, we obtain
\begin{align*}
T(q^2)=\dfrac{T^2(q^4)T(q^8)-2qT(q)T(q^4)T(q^{16})+q^4T(q^4)T^2(q^{16})}{T^2(q)T(q^4)T(q^8)-2q^3T(q)T(q^4)T(q^8)T(q^{16})+q^4T^2(q)T^2(q^{16})}.
\end{align*}
Applying \eqref{T(q)} in the above, we readily arrive at \eqref{r(q)1,2,4,8,16}.

\vspace{.3cm}
\noindent\emph{Proof of \eqref{r(q)1,2,3,6a}}.
From \eqref{gh r(q)1,3 by 2,6 step 2} and \eqref{G^3(q)H(q^3)}, we have
\begin{align*}
&\dfrac{G^3(q)H^3(q)G(q^3)}{G^3(q^2)}-\dfrac{G^2(q^3)H(q^3)}{G(q^6)}\\
&=3q\dfrac{f(-q^2)f^3(-q^{15})}{f^2(-q)f(-q^3)f(-q^{10})}\dfrac{H(q)}{G(q^2)}\\
&=3q\dfrac{f^3(-q^{15})}{f(-q)f(-q^3)f(-q^{5})}\dfrac{G(q)H^2(q)}{G^2(q^2)H(q^2)}\\
&=\dfrac{G(q)H^2(q)}{G^2(q^2)H(q^2)}\bigg(G^3(q)H(q^3)-H^3(q)G(q^3)\bigg),
\end{align*}
which on simplification gives
\begin{align}
\label{G^3(q)H^3(q)H(q^2)G(q^3)}&G^3(q)H^3(q)H(q^2)G(q^3)-G^4(q)H^2(q)G(q^2)H(q^3)+G(q)H^5(q)G(q^2)G(q^3)\\&=\dfrac{G^3(q^2)H(q^2)G^2(q^3)H(q^3)}{G(q^6)}.\notag
\end{align}

Similarly, from \eqref{gh r(q)1,3 by 2,6 step 3} and \eqref{G^3(q)H(q^3)}, we find that
\begin{align}
\label{G^3(q)H^3(q)G(q^2)H(q^3)}&G^3(q)H^3(q)G(q^2)H(q^3)+G^2(q)H^4(q)H(q^2)G(q^3)-G^5(q)H(q)H(q^2)H(q^3)\\&=\dfrac{G(q^2)H^3(q^2)G(q^3)H^2(q^3)}{H(q^6)}.\notag
\end{align}

Dividing \eqref{G^3(q)H^3(q)H(q^2)G(q^3)} by \eqref{G^3(q)H^3(q)G(q^2)H(q^3)}, we have
\begin{align}
\label{G^3(q)H^3(q)H(q^2)G(q^3)a}&\dfrac{G^3(q)H^3(q)H(q^2)G(q^3)-G^4(q)H^2(q)G(q^2)H(q^3)+G(q)H^5(q)G(q^2)G(q^3)}{G^3(q)H^3(q)G(q^2)H(q^3)+G^2(q)H^4(q)H(q^2)G(q^3)-G^5(q)H(q)H(q^2)H(q^3)}\\&=\dfrac{G^2(q^2)G(q^3)H(q^6)}{H^2(q^2)H(q^3)G(q^6)}.\notag
\end{align}
Dividing the numerator and denominator of the left-hand side of the above identity by $G^4(q)H^2(q)G(q^2)H(q^3)$ and then employing \eqref{R-GH}, we obtain 
\begin{align*}
\dfrac{\dfrac{R(q)R(q^2)}{R(q^3)}-1+\dfrac{R^3(q)}{R(q^3)}}{R(q)+\dfrac{R^2(q)R(q^2)}{R(q^3)}-\dfrac{R(q^2)}{R(q)}}=\dfrac{R(q^6)}{R^2(q^2)R(q^3)},
\end{align*}
which can be easily reduced to \eqref{r(q)1,2,3,6a}.

\vspace{.3cm} \noindent\emph{Proof of \eqref{r(q)1,2,3,6,12}}.
Comparing \eqref{n1d1} and \eqref{n2d2}, we have
\begin{align}\label{r123612 step 1}
&\dfrac{H(q)H^2(q^3)G(q^{12})+q G(q)G^2(q^3)H(q^{12})}{G(q)G(q^3)H(q^6)G(q^{12})+q^2H(q)H(q^3)G(q^6)H(q^{12})}\\
&\quad=\dfrac{UG(q^2)H(q^3)G(q^{12})-
q VH(q^2)G(q^3)H(q^{12})}{VH(q^2)H(q^6)G(q^{12})
-q^2 UG(q^2)G(q^6)H(q^{12})}.\notag
\end{align}

Now, we rewrite the terms on the numerators and denominators of the above as follows:
\begin{align*}
&H(q)H^2(q^3)G(q^{12})+q G(q)G^2(q^3)H(q^{12})= G(q)G^2(q^3)H(q^{12})\bigg(\dfrac{T(q)T^2(q^3)}{T(q^{12})}+q\bigg),\\
&G(q)G(q^3)H(q^6)G(q^{12})+q^2H(q)H(q^3)G(q^6)H(q^{12})\\
&\quad= H(q)H(q^3)G(q^6)H(q^{12})\bigg(\dfrac{T(q^6)}{T(q)T(q^3)T(q^{12})}+q^2\bigg),\\
&UG(q^2)H(q^3)G(q^{12})-
q VH(q^2)G(q^3)H(q^{12})\\
&=\bigg( H(q^2)G(q^3)H(q^{12})\bigg)\bigg( G^3(q)H^3(q)H(q^6)-H^3(q^2)G(q^3)H(q^3)\bigg)\\
&\quad\times\bigg(\dfrac{T(q^3)}{T(q^2)T(q^{12})}-
q 
\dfrac{G^3(q)H^3(q)G(q^6)-G^3(q^2)G(q^3)H(q^3)}{G^3(q)H^3(q)H(q^6)-H^3(q^2)G(q^3)H(q^3)} \bigg),\\
&VH(q^2)H(q^6)G(q^{12})
-q^2 UG(q^2)G(q^6)H(q^{12})\\
&=\bigg(G(q^2)G(q^6)H(q^{12})\bigg)\bigg(G^3(q)H^3(q)H(q^6)-H^3(q^2)G(q^3)H(q^3)\bigg)\\&\quad\times\bigg(\dfrac{T(q^2)T(q^6)}{T(q^{12})}\cdot
\dfrac{G^3(q)H^3(q)G(q^6)-G^3(q^2)G(q^3)H(q^3)}{G^3(q)H^3(q)H(q^6)-H^3(q^2)G(q^3)H(q^3)}-
q^2 \bigg),
\end{align*}
where $T(q)$ is as defined in \eqref{T(q)}.

Using the above expressions in \eqref{r123612 step 1}
and then employing \eqref{gh r(q)1,3 by 2,6}, we find that
\begin{align*}
&\dfrac{G(q)G^2(q^3)}{H(q)H(q^3)G(q^6)}\times\dfrac{\dfrac{T(q)T^2(q^3)}{T(q^{12})}+q}{\dfrac{T(q^6)}{T(q)T(q^3)T(q^{12})}+q^2}\\
&=\quad \dfrac{H(q^2)G(q^3)}{G(q^2)G(q^6)}\times \dfrac{\dfrac{T(q^3)}{T(q^2)T(q^{12})}-q\dfrac{T(q)T(q^3)}{T^2(q^2)T(q^6)}}{\dfrac{T(q^2)T(q^6)}{T(q^{12})}\times\dfrac{T(q)T(q^3)}{T^2(q^2)T(q^6)}-q^2},
\end{align*}
which, by \eqref{T(q)}, yields
\begin{align*}
\dfrac{R(q)R^2(q^3)+R(q^{12})}{R(q^6)+R(q)R(q^3)R(q^{12})}=\dfrac{R(q^3)}{R(q^6)}\times\dfrac{R(q^2)R(q^6)-R(q)R(q^{12})}{R(q)R(q^3)-R(q^2)R(q^{12})}.
\end{align*}
Simplifying the previous identity for $R(q^2)$, we arrive at \eqref{r(q)1,2,3,6,12}.

\vspace{.3cm}
\noindent\emph{Proof of \eqref{r(q)1,2,3,4,6,12,24}}.
Comparing \eqref{n3d3} and \eqref{n4d4}, we have
\begin{align}\label{r123461224 step 1}
&\dfrac{WG(q)H(q^{12})G(q^{24})-q^2 X H(q)G(q^{12})H(q^{24})}{q^4WG(q)G(q^6)H(q^{24})
-X H(q)H(q^6)G(q^{24})}\\
&\quad=\dfrac{SYH(q^{12})G(q^{24})-q^2TZG(q^{12})H(q^{24})}{TZH(q^6)G(q^{24})-q^4SYG(q^6)H(q^{24})}.\notag
\end{align}
We rewrite the terms on the numerators and denominators of the above as follows:
\begin{align*}
&WG(q)H(q^{12})G(q^{24})-q^2 X H(q)G(q^{12})H(q^{24})\\
&=\bigg(G(q)H(q^{12})G(q^{24})\bigg)\bigg(G^3(q)G(q^6)H(q^6)-G^3(q^2)H^3(q^2)G(q^3) \bigg)\\
&\quad\times\bigg(\dfrac{H^3(q)G(q^6)H(q^6)-G^3(q^2)H^3(q^2)H(q^3)}{G^3(q)G(q^6)H(q^6)-G^3(q^2)H^3(q^2)G(q^3)}-q^2\dfrac{T(q)T(q^{24})}{T(q^{12})}\bigg),\\
&q^4WG(q)G(q^6)H(q^{24})
-X H(q)H(q^6)G(q^{24})\\
&=\bigg(G(q)G(q^6)H(q^{24})\bigg)\bigg(G^3(q)G(q^6)H(q^6)-G^3(q^2)H^3(q^2)G(q^3) \bigg)\\
&\quad\times\bigg(q^4\dfrac{H^3(q)G(q^6)H(q^6)-G^3(q^2)H^3(q^2)H(q^3)}{G^3(q)G(q^6)H(q^6)-G^3(q^2)H^3(q^2)G(q^3)}-
 \dfrac{T(q)T(q^6)}{T(q^{24})}\bigg),\\
&SYH(q^{12})G(q^{24})-q^2TZG(q^{12})H(q^{24})\\
&=\bigg(G^2(q)H(q^2)G(q^3)H^2(q^4)H^2(q^{12})G(q^{24})\bigg)\\&\quad\times\bigg(H^3(q^2)G(q^3)H(q^{12})-G^3(q)H^3(q^4)H(q^6) \bigg)\\
&\quad\times\bigg(\dfrac{G^3(q^2)H(q^3)G(q^{12})-H^3(q)G^3(q^4)G(q^6)}{H^3(q^2)G(q^3)H(q^{12})-G^3(q)H^3(q^4)H(q^6)}
-q^2\dfrac{T^2(q)T(q^3)T(q^{24})}{T(q^2)T^2(q^4)T^2(q^{12})}\bigg),\end{align*}
 \begin{align*}
&TZH(q^6)G(q^{24})-q^4SYG(q^6)H(q^{24})\\
&=\bigg(H^2(q)G(q^2)H(q^3)G^2(q^4)H(q^6)G(q^{12})G(q^{24})\bigg)\\&\quad\times\bigg(H^3(q^2)G(q^3)H(q^{12})-G^3(q)H^3(q^4)H(q^6)\bigg)\\
&\quad\times\bigg(1-q^4\dfrac{T(q^2)T^2(q^4)T(q^{12})T(q^{24})}{T^2(q)T(q^3)T(q^6)}\cdot\dfrac{G^3(q^2)H(q^3)G(q^{12})-H^3(q)G^3(q^4)G(q^6)}{H^3(q^2)G(q^3)H(q^{12})-G^3(q)H^3(q^4)H(q^6)}\bigg).
\end{align*}

Using the above expressions in \eqref{r123461224 step 1}
and then employing \eqref{gh r(q)1,2,3,6} and \eqref{gh r(q)1 by 2,4}, we obtain
\begin{align*}
&\dfrac{H(q^{12})G(q^{24})}{G(q^6)H(q^{24})}\cdot\dfrac{q^2 T^2(q)T(q^2)T(q^3)T(q^6)-q^2\dfrac{T(q)T(q^{24})}{T(q^{12})}}{q^6 T^2(q)T(q^2)T(q^3)T(q^6)-\dfrac{T(q)T(q^6)}{T(q^{24})}}\\
&\quad=\dfrac{G^2(q)H(q^2)G(q^3)H^2(q^4)H^2(q^{12})}{H^2(q)G(q^2)H(q^3)G^2(q^4)H(q^6)G(q^{12})}\\&\quad\times\dfrac{q^2\dfrac{T(q)}{T(q^2)T(q^4)}-q^2\dfrac{T^2(q)T(q^3)T(q^{24})}{T(q^2)T^2(q^4)T^2(q^{12})}}{1-q^6\dfrac{T(q^2)T^2(q^4)T(q^{12})T(q^{24})}{T^2(q)T(q^3)T(q^6)}\cdot\dfrac{T(q)}{T(q^2)T(q^4)}},
\end{align*}
which, by \eqref{T(q)}, yields 
\begin{align*}
\dfrac{R(q)R(q^2)R(q^3)R(q^6)R(q^{12})-R(q^{24})}{R(q)R(q^2)R(q^3)R(q^6)R(q^{24})-R(q^6)}=\dfrac{R(q^4)R^2(q^{12})-R(q)R(q^3)R(q^{24})}{R(q)R(q^3)R(q^{6})-R(q^4)R(q^{12})R(q^{24})}.
\end{align*}
Simplifying the above for $R(q^2)$, we readily arrive at \eqref{r(q)1,2,3,4,6,12,24}.

\section{Partition-theoretic identities}\label{sec7}

In this section, we present some simple partition-theoretic results arising from the identities in Theorems \ref{GH-new} and \ref{GH-new-ND-R}. 

We first note that $$\dfrac{1}{(q^r;q^s)_\infty^k }$$ is the generating function for the number of partitions of a positive integer into parts that are congruent to $r~(\textup{mod}~s)$ and having $k$ colors. 

For brevity, we define 
\begin{align*}(a_1,a_2,\ldots, a_k;q)_\infty:&=(a_1;q)_\infty(a_2;q)_\infty\cdots(a_k;q)_\infty\\\intertext{and for positive integers $r$ and $s$ with $r<s$,}
(q^{r\pm};q^s)_\infty :&=(q^r, q^{s-r};q^s)_\infty.
\end{align*}

The following theorem arises from \eqref{new 1,2,3,6 first}.
\begin{theorem}
Let $p_1(n)$ denote the number of partitions of $n$ into parts congruent to $\pm 1, \pm 5, \pm 11, \pm 12 ~\pmod {30}$, where the
parts congruent to $\pm  5,  \pm 12 ~\pmod {30}$ have two colors. Let $p_2(n)$ denote the number of partitions of $n$ into parts congruent to $\pm 5, \pm 6, \pm 7, \pm 13 ~\pmod {30}$, where the
parts congruent to $\pm  5,  \pm 6 ~\pmod {30}$ have two colors. Let $p_3(n)$ denote the number of partitions of $n$ into parts congruent to $\pm 1, \pm 6, \pm7, \pm 11, \pm 12, \pm 13 ~\pmod {30}$. Then, for each positive integer $n$,
\begin{align}\label{part-int1} p_1(n)-p_2(n)=p_3(n-1).
\end{align}
\end{theorem}

\begin{proof}Applying \eqref{rrident1} and \eqref{rrident2} in \eqref{new 1,2,3,6 first}, we find that
\begin{align*}
\dfrac{(q^{2\pm};q^5)_\infty(q^{6\pm};q^{30})_\infty
}{(q^{2\pm};q^{10})_\infty(q^{3\pm};q^{15})_\infty}-\dfrac{(q^{1\pm};q^5)_\infty(q^{12\pm};q^{30})_\infty
}{(q^{4\pm};q^{10})_\infty(q^{6\pm};q^{15})_\infty}=q(q^{5\pm};q^{30})_\infty^2.
\end{align*}
Dividing both sides of the above by $(q^{1\pm},q^{2\pm};q^5)_\infty(q^{5\pm};q^{30})_\infty^2(q^{6\pm},q^{12\pm} ;q^{30})_\infty$, reducing all the products into base $q^{30}$, and then cancelling the common terms, we obtain
\begin{align*}
&\dfrac{1}{(q^{1\pm},q^{5\pm},q^{5\pm},q^{11\pm},q^{12\pm},q^{12\pm};q^{30})_\infty
}-\dfrac{1}{(q^{5\pm},q^{5\pm},q^{6\pm},q^{6\pm},q^{7\pm},q^{13\pm};q^{30})_\infty}\\
&\quad=\dfrac{q}{(q^{1\pm},q^{6\pm},q^{7\pm},q^{11\pm},q^{12\pm},q^{13\pm};q^{30})_\infty
},
\end{align*}
which is equivalent to 
\begin{align*}
&\sum_{n=0}^\infty p_1(n)q^n-\sum_{n=0}^\infty p_2(n)q^n=\sum_{n=0}^\infty p_3(n)q^{n+1}.
\end{align*}
Equating the coefficients of $q^n$ for $n\ge1$, we readily arrive at \eqref{part-int1} to finish the proof.
\end{proof}

In a similar way, one can easily derive the following two theorems from \eqref{new 1,2,3,6 second} and \eqref{n1d1}, respectively.
\begin{theorem}
Let $p_4(n)$ denote the number of partitions of $n$ into parts congruent to $\pm 2, \pm 3, \pm 7, \pm 8,\pm 12,\pm 13 ~\pmod {30}$, where the
parts congruent to $\pm  2, \pm 3,\pm 8, \pm 12 ~\pmod {30}$ have two colors. Let $p_5(n)$ denote the number of partitions of $n$ into parts congruent to $\pm 1, \pm 4,\pm 6, \pm 9, \pm 11,\pm 14 ~\pmod {30}$, where the
parts congruent to $\pm  4,  \pm 6, \pm 9, \pm 14~\pmod {30}$ have two colors. Let $p_6(n)$ denote the number of partitions of $n$ into parts congruent to $\pm 2, \pm 3,\pm 4$, $\pm 5,\pm 6, \pm8,\pm 9, \pm 12, \pm 14 ~\pmod {30}$ where the parts congruent to $\pm  5~\pmod {30}$ have two colors. Then, for each positive integer $n>1$,
$$p_4(n)-p_5(n-2)=p_6(n).$$ 
\end{theorem}

\begin{theorem}
Let $p_7(n)$ denote the number of partitions of $n$ into parts congruent to $\pm 2, \pm 6, \pm 7, \pm 8,\pm 9,\pm 10,\pm 12,\pm 13,\pm 17,\pm 21,\pm 22,\pm 23,\pm 24,\pm 28 ~\pmod {60}$, where the parts congruent to $\pm  10 ~\pmod {60}$ have two colors. Let $p_8(n)$ denote the number of partitions of $n$ into parts congruent to $\pm 1, \pm 3, \pm 4$, $\pm 10,\pm 11,\pm 12,\pm 14,\pm 16,\pm 18,$ $\pm 19,\pm 24,\pm 26,\pm 27,\pm 29 ~\pmod {60}$, where the
parts congruent to $\pm  10~\pmod {60}$ have two colors.  Let $p_9(n)$ denote the number of partitions of $n$ into parts congruent to $\pm 1, \pm 4, \pm 5, \pm 6,\pm 11,\pm 12,\pm 14,\pm 16$, $\pm 18, \pm 19,\pm 25,\pm 26, \pm 29 ~\pmod {60}$, where the
parts congruent to $\pm  12, \pm 18~\pmod {60}$ have two colors. Let $p_{10}(n)$ denote the number of partitions of $n$ into parts congruent to $\pm 2, \pm 5,  \pm 6,\pm 7,\pm 8,\pm 13,\pm 17,\pm 18,\pm 22,$ $\pm 23,\pm 24,\pm 25,\pm 28 ~\pmod {60}$, where the
parts congruent to $\pm  6, \pm 24~\pmod {60}$ have two colors. Then, for each positive integer $n>1$,
$$p_7(n)+p_8(n-1)=p_9(n)+p_{10}(n-2).$$ 
\end{theorem}

Similar partition-theoretic identities may also be derived from the other identities in Theorem \ref{GH-new-ND-R}. However, we omit those considerably lengthy statements.

\section{Concluding remarks}

In this paper, we find some new modular identities for the Rogers-Ramanujan continued fraction, namely, Thereom \ref{theoremnew1} and Theorem \ref{theoremnew2} by using two approaches.  The identities in Thereom \ref{theoremnew1} and  \eqref{r(q)1,2,3,6}--\eqref{r(q)1,2,4,8,16} in Thereom \ref{theoremnew2} are derived by using dissections of theta functions. Identities \eqref{r(q)1,2,3,6a}--\eqref{r(q)1,2,3,4,6,12,24} in Thereom \ref{theoremnew2} are derived from some  identities for the Rogers-Ramanujan functions arising from the quintuple product identity \eqref{qtpi}. In Lemmas \ref{lem-q78}--\ref{lem-1317}, some theta function identities derived from \eqref{qtpi} are given. We then use those identities to  deduce the identities in Theorem \ref{GH-new-ND-R} involving the Rogers-Ramanujan functions, which then imply the identities \eqref{r(q)1,2,3,6a}--\eqref{r(q)1,2,3,4,6,12,24} for the Rogers-Ramanujan continued fraction in Theorem \ref{GH-new-ND-R}. The next set of theta function identities analogous to those in Lemmas \ref{lem-q78}--\ref{lem-1317} may be explored. However, it would be challenging to derive identities for the Rogers-Ramanujan continued fraction analogous to \eqref{r(q)1,2,3,6a}--\eqref{r(q)1,2,3,4,6,12,24} from them. In the following, we briefly outline a few steps in that direction.

First, setting $q=q^{10}$ and $B=q$, $q^3$, $q^7$ and $q^9$, respectively, we find that
\begin{align}
\label{f(q^13,q^47)}f\left(q^{13},q^{47}\right)-q^2f\left( q^{7},q^{53}\right)&=f(-q^{20})\dfrac{f\left(-q^{2},-q^{18}\right)}{f\left(q^{9},q^{11}\right)},\\
\label{f(q^19,q^41)}f\left(q^{19},q^{41}\right)-q^6f\left( q,q^{59}\right)&=f(-q^{20})\dfrac{f\left(-q^{6},-q^{14}\right)}{f\left(q^{7},q^{13}\right)},\\
\label{f(q^29,q^31)}f\left(q^{29},q^{31}\right)-q^3f\left( q^{11},q^{49}\right)&=f(-q^{20})\dfrac{f\left(-q^{6},-q^{14}\right)}{f\left(q^{3},q^{17}\right)},\\
\label{f(q^23,q^37)}f\left(q^{23},q^{37}\right)-q f\left( q^{17},q^{43}\right)&=f(-q^{20})\dfrac{f\left(-q^{2},-q^{18}\right)}{f\left(q,q^{19}\right)}.
\end{align}
Proceeding as in the proof of Lemma \ref{lem-78}, it can be shown that
\begin{align}
\label{f(q^13,q^47)^3}f^3\left(q^{13},q^{47}\right)-q^6f^3\left( q^{7},q^{53}\right)&=f^3(-q^{20})\dfrac{f\left(-q^{6},-q^{54}\right)}{f\left(q^{27},q^{33}\right)},\\
\label{f(q^19,q^41)^3}f^3\left(q^{19},q^{41}\right)-q^{18}f\left( q,q^{59}\right)&=f^3(-q^{20})\dfrac{f\left(-q^{18},-q^{42}\right)}{f\left(q^{21},q^{39}\right)},\\
\label{f(q^29,q^31)^3}f^3\left(q^{29},q^{31}\right)-q^9f\left( q^{11},q^{49}\right)&=f^3(-q^{20})\dfrac{f\left(-q^{18},-q^{42}\right)}{f\left(q^{9},q^{51}\right)},\\
\label{f(q^23,q^37)^3}f^3\left(q^{23},q^{37}\right)-q^3f^3\left( q^{17},q^{43}\right)&=f^3(-q^{20})\dfrac{f\left(-q^{6},-q^{54}\right)}{f\left(q^3,q^{57}\right)}.
\end{align}
However, from \eqref{f(a,b)+f(-a-b)} and \eqref{f(a,b)-f(-a-b)}, we have
\begin{align*}
f\left(q^{9},q^{11}\right)&=\dfrac{1}{2}\left(f\left(q^{2},q^{3}\right)+f\left(-q^{2},-q^{3}\right)\right),\\
f\left(q^{7},q^{13}\right)&=\dfrac{1}{2}\left(f\left(q,q^{4}\right)+f\left(-q,-q^{4}\right)\right),\\
f\left(q,q^{19}\right)&=\dfrac{1}{2q^2}\left(f\left(q^{2},q^{3}\right)-f\left(-q^{2},-q^{3}\right)\right),\\
f\left(q^{3},q^{17}\right)&=\dfrac{1}{2q}\left(f\left(q,q^{4}\right)-f\left(-q,-q^{4}\right)\right).
\end{align*}
Therefore, \eqref{f(q^13,q^47)}--\eqref{f(q^23,q^37)^3} may be recast as 
\begin{align}
\label{f(q^13,q^47)a}f\left(q^{13},q^{47}\right)-q^2f\left( q^{7},q^{53}\right)&=f(-q^{20})\dfrac{2f\left(-q^{2},-q^{18}\right)}{f\left(q^{2},q^{3}\right)+f\left(-q^{2},-q^{3}\right)},\\
\label{f(q^19,q^41)a}f\left(q^{19},q^{41}\right)-q^6f\left( q,q^{59}\right)&=f(-q^{20})\dfrac{2f\left(-q^{6},-q^{14}\right)}{f\left(q,q^{4}\right)+f\left(-q,-q^{4}\right)},\\
\label{f(q^29,q^31)a}f\left(q^{29},q^{31}\right)-q^3f\left( q^{11},q^{49}\right)&=f(-q^{20})\dfrac{2qf\left(-q^{6},-q^{14}\right)}{f\left(q,q^{4}\right)-f\left(-q,-q^{4}\right)},\\
\label{f(q^23,q^37)a}f\left(q^{23},q^{37}\right)-q f\left( q^{17},q^{43}\right)&=f(-q^{20})\dfrac{2q^2f\left(-q^{2},-q^{18}\right)}{f\left(q^2,q^3\right)-f\left(-q^2,-q^{3}\right)},\\
\label{f(q^13,q^47)^3a}f^3\left(q^{13},q^{47}\right)-q^6f^3\left( q^{7},q^{53}\right)&=f^3(-q^{20})\dfrac{2f\left(-q^{6},-q^{54}\right)}{f\left(q^{6},q^{9}\right)+f\left(-q^{6},-q^{9}\right)},\end{align}
\begin{align}
\label{f(q^19,q^41)^3a}f^3\left(q^{19},q^{41}\right)-q^{18}f\left( q,q^{59}\right)&=f^3(-q^{20})\dfrac{2f\left(-q^{18},-q^{42}\right)}{f\left(q^{3},q^{12}\right)+f\left(-q^{3},-q^{12}\right)},\\
\label{f(q^29,q^31)^3a}f^3\left(q^{29},q^{31}\right)-q^9f\left( q^{11},q^{49}\right)&=f^3(-q^{20})\dfrac{2q^3f\left(-q^{18},-q^{42}\right)}{f\left(q^{3},q^{12}\right)-f\left(-q^{3},-q^{12}\right)},\\
\label{f(q^23,q^37)^3a}f^3\left(q^{23},q^{37}\right)-q^3f^3\left( q^{17},q^{43}\right)&=f^3(-q^{20})\dfrac{2q^6f\left(-q^{6},-q^{54}\right)}{f\left(q^{6},q^{9}\right)-f\left(-q^{6},-q^{9}\right)}.
\end{align}
Then proceeding as in the proofs of \eqref{n2d2} and \eqref{gh r(q)1,3 by 2,6}, we have the following identities from \eqref{f(q^13,q^47)a} and \eqref{f(q^13,q^47)^3a}, which are somewhat analogous to \eqref{n2d2 step 1} and \eqref{gh r(q)1,3 by 2,6 step 2}: 
\begin{align}
\label{rem1}&f^3(-q)f(-q^{60})G^3(q^2)H(q^6)H(q^{12})\left(G(q)H^2(q)+G(q)H(q^2)\right)^3\\
&\quad-4f(-q^3)f^3(-q^{20})H^3(q^2)H^3(q^4)G(q^6)\left(G(q^3)H^2(q^3)+G(q^3)H(q^6)\right)\notag\\
&=\dfrac{3q^2f^2(-q)f(-q^{3})}{f(-q^{20})}G^2(q^2)H(q^2)H(q^4)G(q^6)\left(G(q)H^2(q)+G(q)H(q^2)\right)^2\notag\\
&\quad\times\left(G(q^3)H^2(q^3)+G(q^3)H(q^6)\right)f(q^{13},q^{47})f(q^7,q^{53}),\notag\\
\label{rem2}& f^3(-q)f(-q^{60})G^3(q^2)H(q^6)H(q^{12})\left(G(q)H^2(q)+G(q)H(q^2)\right)^3\\
&\quad-4f(-q^3)f^3(-q^{20})H^3(q^2)H^3(q^4)G(q^6)\left(G(q^3)H^2(q^3)+G(q^3)H(q^6)\right)\notag\\
&=\dfrac{3q^2f^3(-q)f(-q^{10})f(-q^{30})f^3(-q^{60})}{f(-q^2)f(-q^6)f^2(-q^{20})}H(q^4)\left(G(q)H^2(q)+G(q)H(q^2)\right)^2\notag\\
&\quad\times\left(G^2(q)H(q)+G(q^2)H(q)\right).\notag
\end{align}
In a similar fashion, using  \eqref{f(q^19,q^41)a} and \eqref{f(q^19,q^41)^3a}, it can be shown that
\begin{align}
\label{rem3}&f^3(-q)f(-q^{60})H^3(q^2)G(q^6)G(q^{12})\left(G^2(q)H(q)+H(q)G(q^2)\right)^3\\
&\quad-4f(-q^3)f^3(-q^{20})G^3(q^2)G^3(q^4)H(q^6)\left(G^2(q^3)H(q^3)+H(q^3)G(q^6)\right)\notag\\
&=\dfrac{3q^6f^2(-q)f(-q^{3})}{f(-q^{20})}G(q^2)H^2(q^2)G(q^4)H(q^6)\left(G^2(q)H(q)+H(q)G(q^2)\right)^2\notag\\
&\quad\times\left(G^2(q^3)H(q^3)+H(q^3)G(q^6)\right)f(q^{19},q^{41})f(q,q^{59})\notag\\
\intertext{and}
\label{rem4}& f^3(-q)f(-q^{60})H^3(q^2)G(q^6)G(q^{12})\left(G^2(q)H(q)+H(q)G(q^2)\right)^3\\
&\quad-4f(-q^3)f^3(-q^{20})G^3(q^2)G^3(q^4)H(q^6)\left(G^2(q^3)H(q^3)+H(q^3)G(q^6)\right)\notag\\
&=\dfrac{3q^4f^3(-q)f(-q^{10})f(-q^{30})f^3(-q^{60})}{f(-q^2)f(-q^6)f^2(-q^{20})}G(q^4)
\left(G^2(q)H(q)+H(q)G(q^2)\right)^2\notag\\
&\quad\times\left(G(q)H^2(q)-G(q)H(q^2)\right).\notag
\end{align}
Similar identities may be derived from \eqref{f(q^29,q^31)a},\eqref{f(q^23,q^37)a}, \eqref{f(q^29,q^31)^3a}, and \eqref{f(q^23,q^37)^3a}. However, the appearance of terms, like  $f(q^{13},q^{47})f(q^7,q^{53})$ in \eqref{rem1} and $f(q^{19},q^{41})f(q,q^{59})$ in \eqref{rem3}, might be a hindrance in deriving identities analogous to \eqref{r(q)1,2,3,6a}--\eqref{r(q)1,2,3,4,6,12,24}.

\section*{Acknowledgement}
 %The authors are thankful to the anonymous referee for his/her helpful comments on the paper. 
 The second author was partially supported by Council of Scientific \& Industrial Research (CSIR), Government of India under CSIR-JRF scheme (Grant No. 09/0796(12991)/2021-EMR-I). The author thanks the funding agency.

\end{document}